\documentclass[10pt,a4paper]{amsart}

\usepackage{setspace}
\usepackage[utf8x]{inputenc}
\usepackage[T1]{fontenc}

\usepackage{helvet}

\usepackage{amsthm}
\usepackage{amsmath}
\usepackage{verbatim}
\usepackage{amssymb} 
\usepackage{resizegather}
\usepackage{tikz}
\usepackage{tikz-cd}
\usepackage{babel} 
\usepackage{mathrsfs}
\usepackage{graphicx}
\usepackage{epigraph}

\usepackage{pgfplots}
\usepgfplotslibrary{fillbetween}
\pgfplotsset{compat=1.13}

\usepackage[bottom]{footmisc}

\usepackage{hyperref}

\newtheorem{theorem}{Theorem}[section]

\newtheorem{proposition}{Proposition}[theorem]

\newtheorem{lemma}{Lemma}[theorem]

\newtheorem{corollary}{Corollary}[theorem]

\newtheorem{remark}{Remark}
\newtheorem{example}{Example}
\newtheorem{definition}{Definition}

\title{On Isolated Real Singularities I}
\author{Lars Andersen}

\begin{document}

\maketitle

\begin{abstract} 
This article and its successor concern the topology of real isolated hypersurface singularities. We prove in Theorem \hyperref[main theorem one]{\ref*{main theorem one}} that after attaching a certain number of handles the real Milnor fibres become contractible, with each handle corresponding to a critical point of a $\mathbb{R}$-morsification (Definition \hyperref[definition morsification]{\ref*{definition morsification}}); in particular one recovers the formula of Khimshiashvili (\cite[Theorem 2.3]{Khimshiashvili}) for the Euler characteristic of the Milnor fibres.\\
We then give sufficient conditions, in Theorem \hyperref[main theorem two]{\ref*{main theorem two}} for having that the integer homology groups of the real Milnor fibres are isomorphic to the homology groups of a bouquet of spheres. It follows that the homology groups of the real Milnor fibres are uniquely determined by the real vanishing cycles (Definition \hyperref[vanishing cycle]{\ref*{vanishing cycle}}) in this case.
\end{abstract}
\begin{center}
    \text{Classification: }\textbf{53-XX, 57-XX}
\end{center}

\section{Introduction}

If one of the principal objectives in the study of the singularities of maps of algebraic varieties is to obtain reasonable classifications of these, then an adequate understanding of calculable invariants of topological or algebraic nature is desirable first of all, and secondly, suitable equivalence relations preserving a given such invariant. For singularities defined over the field of complex numbers several such invariants has been used for the classification problem, amongst others the Milnor number and the Tjurina number, and several equivalence relations, amongst others the contact equivalence and the $\mu$-constant equivalence relations, preserving one or more of these invariants. The situation is as follows. Let $f: (\mathbb{R}^{n+k},0)\to(\mathbb{R}^k,0)$ be a germ of real analytic function, choose a representative $f: M\to\mathbb{R}^k$ on a small neighborhood $M\subset\mathbb{R}^{n+k}$ containing the origin and let $\text{Jac}(f)$ denote the Jacobian matrix of this representative. Then

\begin{theorem}[{\cite[Theorem 4.2]{Seade}}]\label{milnor one} If the Jacobian matrix $\text{Jac}(f)$ has rank $k$ in any $x\neq 0$, then there exists $\delta_0>0$ such that for any $\delta\in (0,\delta_0]$ there exists $\epsilon_0>0$ such that for any $\epsilon\in (0,\epsilon_0]$ the restriction
\begin{equation}\label{projection}
f: f^{-1}(\mathbb{B}_{\epsilon}^{k}\setminus\{0\})\cap \mathbb{B}_{\delta}^{n+k}\to \mathbb{B}_{\epsilon}^{k}\setminus\{0\}
\end{equation}
is a locally trivial topological fibration where $\mathbb{B}_{\delta}^{n+k}\subset \mathbb{R}^{n+k}$ denotes the open $(n+k)$-dimensional ball centered at the origin and of radius $\delta$ and where $\mathbb{B}_{\epsilon}^{k}\subset\mathbb{R}^k$ denotes the open $k$-dimensional ball centered at the origin and of radius $\epsilon$.
\end{theorem}

In the case of hypersurface singularities (i.e $k=1$) there are therefore up to\\
diffeomorphism two fibres 
$$\mathcal{F}_{\eta}^{\pm}:= f^{-1}(\pm\eta)\cap \mathbb{B}_{\delta}^{n+1},\qquad 0<\eta\leq\epsilon$$
of the locally trivial fibration \hyperref[projection]{\ref*{projection}}. These are referred to as the \emph{positive} and \emph{negative} open \emph{Milnor fibres} of $f$ at the origin. The closures $\bar{\mathcal{F}}_{\eta}^{\pm}$ for the euclidean topology of the positive and negative open Milnor fibres are then manifolds with boundaries and are referred to as the positive and negative \emph{closed Milnor fibres} of $f$ at the origin.\\
It is a delicate problem to effectively calculate topological invariants for the Milnor fibres of real singularities\footnote{More generally, in the words of Réné Thom as quoted in the introduction (p. 5) to the book \cite{roy}, real algebraic geometry, in Thom's view an example of mathematics pertaining to phenomena where ``reality plays an essential role'', has historically been neglected in lieu of the ``even too beautiful'' complex geometry. This might arguably explain at least in part why relatively little is known of the topology of real singularities.} and as a consequence, the different classification problems at hand, are mostly completely open at the present date. To approach this problem we shall by and large to adapt the proof of E. Brieskorn \cite[Appendix]{brieskorn} of Milnors theorem \cite[Theorem 7.2]{Milnor} by replacing arguments which fails in the real case by arguments coming from Stratified Morse Theory \cite{Goresky} as constructed by Goresky and MacPherson. 
\subsection*{Acknowledgments} The author wish to express his sincere gratitude to the\\
Laboratoire de Mathématiques at USMB and especially to his thesis supervisors Georges Comte and Michel Raibaut. The material formed in this article and its successor formed a part of the authors thesis. On behalf of all authors, the corresponding author states that there is no conflict of interest.
\subsection*{Data Sharing} Data sharing not applicable to this article as no datasets were generated or analysed during the current study.
\section{Isolated Real Hypersurface Singularities}
\subsection{$\mathbb{R}$-Morsifications}
A principal idea which presents itself in geometry is that of deforming a given object as to eliminate nongeneric properties of the object in question whilst conserving a given geometrical invariant. One example is the moving lemma in intersection theory whereby one can perturb two given cycles of algebraic varieties as to intersect transversally, whilst preserving their intersection number. Another example presents itself in the theory of singularities, where one can deform a germ of holomorphic functions, defining an isolated hypersurface singularity, as to split the singularity into isolated quadratic singularities\footnote{This preserves the \emph{total} Milnor number}.\\ 

We give the following definition.\\

\begin{definition}\label{definition morsification} Let $f:(\mathbb{R}^{n+1},0)\to(\mathbb{R},0)$ be a germ of polynomial maps. A $\mathbb{R}$-morsification of $f$ is a representative 
$$F: M\times U\to\mathbb{R}$$ 
of a polynomial map germ $F: (\mathbb{R}^{n+1}\times \bar{\mathcal{I}}, 0)\to(\mathbb{R},0)$ such that
\begin{enumerate}
\item $\bar{\mathcal{I}}\cong [-1,1]^k$ for some $k\in\mathbb{N}$,
\item $F(x,0)=f(x)$ as germs of polynomial maps,
\item There exists a dense subset $V\subset U\setminus\{0\}$ such that for all $t\in V$, the function $f_t:M\to\mathbb{R}$ is a Morse function with distinct critical values.
\end{enumerate}
An $\mathbb{R}$-morsification $F$ is \emph{strong} if for any $t\in V$ the function $f_t: M\to\mathbb{R}$ has $\mu(f_{\mathbb{C}},0)$ critical points and it is \emph{weak} if for any $t\in V$ it has strictly less than $\mu(f_{\mathbb{C}},0)$ critical points.
\end{definition}

In the case where the polynomial map germ has an isolated critical point in the origin\footnote{that is: for any representative $f: M\to\mathbb{R}$ the origin is an isolated critical point} an $\mathbb{R}$-morsification of $f$ always exist:

\begin{proposition}\label{simple way} If $f: (\mathbb{R}^{n+1},0)\to(\mathbb{R},0)$ has an isolated critical point in the origin then $f$ has an $\mathbb{R}$-morsification.
\end{proposition}
\begin{proof} The proof of \cite[Proposition 3.18]{Ebeling} applies almost verbatim, but since it is constructive we recall the argument. Namely let $f: M\to\mathbb{R}$ with $f(0)=0$ be a representative of the germ $f: (\mathbb{R}^{n+1},0)\to(\mathbb{R},0)$, for some neighborhood $M\subset\mathbb{R}^{n+1}$ of the origin. Consider the real algebraic subset $\mathcal{D}\subset\mathbb{R}^{n+1}$ of critical values of the gradient  
$$\text{grad}(f): \mathbb{R}^{n+1}\to\mathbb{R}^{n+1}.$$
Its complement is a dense subset by the Brown-Sard Theorem \cite{Kosinski}{Theorem 3.1}. We claim that the tangent cone $T\mathcal{D}\subset\mathbb{R}^{n+1}$ of $\mathcal{D}$ satisfies $T\mathcal{D}\neq \mathbb{R}^{n+1}$. Indeed suppose not. Then the second Fréchet derivative $d(\text{grad}(f))$ has constant rank $k\leq n+1$ in a neighborhood $V\subset M$ of the origin. By the Constant Rank Theorem \cite[Theorem 1.2]{Kosinski}, 
$$\text{grad}(f)(x_1,\dots, x_{n+1})=(x_1,\dots, x_k, 0,\dots, 0)$$
in local coordinates. If $k=n+1$ then $\text{grad}(f)$ is the identity map and if $k\leq n$ then 
$$\partial f/\partial x_i=0,\qquad i=k+1,\dots, n+1$$
hence $f$ is constant; in either case one gets a contradiction to the fact that the origin was an isolated critical point. Therefore $T\mathcal{D}$ has positive codimension in $\mathbb{R}^{n+1}$ which implies that we can choose a union 
$$l=\{(ta_1,\dots, ta_{n+1})\ |\ t\in (-a, 0)\cup (0, a)\},\qquad a\in\mathbb{R}$$
of two line segments such that each line segment does not meet any critical value. Letting  
$$f_t=f+t\sum_{i=1}^{n+1} a_i x_i$$
then the critical points of $f_t$ are exactly the points $p\in M$ for which 
$$\text{grad}(f)(p)=(t a_1,\dots, t a_{n+1}).$$
This is a noncritical value of the gradient $\text{grad}(f)$ hence the Hessian matrix $\text{Hess}(f_t)(p)$ at $p$ is invertible. But this means that $p$ is a nondegenerate critical point of $f_t$. Therefore if one defines the subset $V\subset\mathbb{R}$ by $V=(-a, 0)\cup (0,a)$ then for any $t\in V$ one has that $f_t$ has only nondegenerate critical points in $M$. Up to changing the vector $(a_1,\dots, a_{n+1})$ by a small perturbation one can assume that $f_t$ has only nondegenerate critical points with distinct critical values for any $t\in V$. Then the representative 
$$F: M\times (-a, a)\to\mathbb{R},\qquad F(x,t)=f_t(x)$$ 
of the germ 
$$F: (\mathbb{R}^{n+1}\times [-a, a], 0)\to(\mathbb{R},0)$$
satisfies the conditions of the Definition \hyperref[definition morsification]{\ref*{definition morsification}} of a $\mathbb{R}$-morsification, with $U=(-a, a)$ and $V=U\setminus\{0\}$.
\end{proof}

\begin{remark} As the proof shows one can always obtain a weak $\mathbb{R}$-morsification of a germ of isolated singularity $f: (\mathbb{R}^{n+1}, 0)\to(\mathbb{R},0)$ having the property that $V=U\setminus\{0\}$. Note however that the \emph{number} of critical points of $f_t: M\to\mathbb{R}$ might be different for different choices of $t\in V$. 
\end{remark}

We give a few examples.

\begin{example} Let $f:\mathbb{R}\to\mathbb{R}$ be given by $f(x)=x^4$ and let $F: \mathbb{R}\times\mathbb{R}\to\mathbb{R}$ be $F(x,t)=x^4-tx$. Then
$$f_t'(x)=4x^3-t=0\qquad\Leftrightarrow x^3=t/4$$
and as there are only one real third root of unity there is only one critical point of $f_t:\mathbb{R}\to\mathbb{R}$ for any $t\in\mathbb{R}$. But the Milnor number of $f$ in the origin is $\mu(f,0)=3$. Therefore there is no closed interval $\bar{\mathcal{I}}$ containing the origin such that $F_{|\mathbb{R}\times\bar{\mathcal{I}}}$ is a strong $\mathbb{R}$-morsification.
\end{example}

In this example one can find another one-parameter family which is a strong $\mathbb{R}$-morsification.

\begin{example} Let $F: \mathbb{R}\times I\to\mathbb{R}$ with $I=[0,32/27]$ be given by $F(x,t)=x^4-4tx^2+4t^2x$. Then for any $t\in[0, 32/27)$, the derivative $f_t'=4(x^3-2tx+t^2)$ has discriminant
$$\Delta(f_t')=t^3(32-27t)$$
which is positive when $0<t<32/27$ whence it follows that $f_t$ have exactly $\mu(f,0)=3$ real roots. The set of $t\in (0,32/27]$ such that the critical points of $f_t$ are not Morse is the image $\pi(X)$ of  
$$X=\{(x,t)\ |\ 3x^2-2t=0,\quad x^3-2tx+t^2=0\}\subset \mathbb{R}\times (0,32/27]$$
under the projection $\pi: \mathbb{R}\times\mathbb{R}\to\mathbb{R}$. As $X$ is a proper intersection of two curves it follows that $\dim X=0$ whence $\pi(X)$ has finite cardinality. Therefore $V=(0, 32/27]\setminus \pi(X)$ is dense in $I$ and for any $t\in V$ the functions $f_t$ are Morse. Hence $F$ is a strong $\mathbb{R}$-morsification.
\end{example}

As for the existence of strong $\mathbb{R}$-morsifications let us recall the notion of $M$-morsification, due to Vladimir I. Arnold. 

\begin{definition}[{\cite{Guryunov}}] A $M$-morsification of a germ of holomorphic functions $f: (\mathbb{C}^{n+1},0)\to (\mathbb{C},0)$ with an isolated singularity in the origin and of Milnor number $\mu=\mu(f,0)$, is a representative 
$$F: M\times U\to\mathbb{C}$$
of a germ of holomorphic functions $F: (\mathbb{C}^{n+1} \times \mathbb{C}^{\mu-1},0) \to (\mathbb{C},0)$ such that
\begin{enumerate} 
\item $F(\cdot,0)=f(\cdot)$ as germs of holomorphic functions, 
\item $U\subset \mathbb{R}^{\mu-1}$ and there exists a subset $V\subset U\setminus\{0\}$ of complement of Lebesgue measure zero such that for all $t\in U$ the function $f_t=F(\cdot, t)$ is Morse with exactly $\mu$ real critical points.
\end{enumerate}
A $M$-morsification is said to be \emph{very nice} if moreover the critical values of $f_t$ are distinct, for all $t\in U$.
\end{definition}

Stated differently, an $M$-morsification is defined from a representative of the truncated miniversal deformation\footnote{often called the truncated universal unfolding} of $f$ by restricting the base to the (open) set of real parameters for which all nondegenerate critical points of $f_t$ are real.\\

Arnold established in his article \cite{Arnold} from 1991 the existence of very nice $M$-morsifications for simple singularities of class $A_k$; more precisely he calculated therein the number of connected components  of the space of such morsifications. In 1993 M.R. Entov established the existence of very nice $M$-morsifications for $D_k$-singularities. And in his article \cite{Guryunov} V.V. Guryunov proved their existence for simple singularities of class $E_6, E_7$ and $E_8$, showing that the number of connected components of the space of very nice $M$-morsifications are $82, 768$ and $4056$, respectively. 

\begin{proposition}\label{ADE one} Let $f: (\mathbb{R}^{n+1},0)\to(\mathbb{R},0)$ be a germ of polynomial maps defining an isolated singularity in the origin. Suppose that for any representative $f: M\to\mathbb{R}$, where $M$ is an open neighbourhood of the origin, the extension of scalars $f_{\mathbb{C}}: M_{\mathbb{C}}\to\mathbb{C}$ is a simple isolated singularity of class $ADE$. Then there exists a strong $\mathbb{R}$-morsification of $f$.
\end{proposition}
\begin{proof} 
\begin{enumerate}
\begin{item} Suppose $n=1$. As $f_{\mathbb{C}}: M_{\mathbb{C}}\to\mathbb{C}$ is a simple $ADE$-singularity, by the results \cite{Arnold}, \cite{Entov} and \cite{Guryunov}, there exists a very nice $M$-morsification 
$$F_{\mathbb{C}}: M_{\mathbb{C}}\times U \to\mathbb{C}$$
of $f_{\mathbb{C}}$ where $U\subset\mathbb{R}^{\mu-1}$ is an open neighbourhood of the origin. By definition there exists a subset $V\subset U\setminus\{0\}$ of complement of measure zero such that the functions
$$f_{\mathbb{C},t}:M_{\mathbb{C}}\to\mathbb{C},\qquad t\in V$$
are Morse with $\mu$ real critical points $p_1,\dots, p_{\mu}$ and distinct critical values. There might be other critical points of $f_{\mathbb{C},t}$ which are not real, but as we show now this is not the case. Indeed the complement of a set $A\subset\mathbb{R}^k$ of measure zero is dense, so $V$ is dense. And as $F_{\mathbb{C}}$ is a representative of a miniversal deformation it follows by \cite[Lemma 3.9]{Ebeling} that it is a strict morsification. By \cite[Proposition 3.19]{Ebeling} as $F_{\mathbb{C}}$ is a morsification of $f_{\mathbb{C}}$ it follows that $f_{\mathbb{C},t}$ has all in all $\mu$ nondegenerate critical points. Therefore all the critical points are real. Consequently if one considers the restriction
$$F: M\times U \to\mathbb{R}$$
to the open neighborhood $M$ of the origin in $\mathbb{R}^{2}$ then $F(x,0)=f(x)$ for all $x\in M$ and moreover, for any $t\in V$, the function $f_t=F(\cdot,t)$ is Morse with distinct critical values. Hence $F$ is a strong $\mathbb{R}$-morsification.

\item If $n\geq 1$ and if $f: M\to\mathbb{R}$ is a representative of the germ $f: (\mathbb{R}^{n+1}, 0)\to (\mathbb{R}, 0)$ then one can by assumption write
$$f(x_1,\dots,x_{n+1})=\tilde{f}(x_1,x_2)+\sum_{i=3}^{n+1} \epsilon_i x_i^2,\qquad \epsilon_i\in\{\pm 1\}$$
where the extension of scalars of $\tilde{f}:M\cap \mathbb{R}^2\to\mathbb{R}$ is a simple curve singularity of class $ADE$. Notice that 
$$\mu=\mu(f_{\mathbb{C}},0)=\mu(\tilde{f}_{\mathbb{C}},0).$$
By the step above there exists a strong $\mathbb{R}$-morsification 
$$\tilde{F}: M\cap\mathbb{R}^2 \times U\to\mathbb{R}$$
of $\tilde{f}$, where $U\subset [-\delta, \delta]^{\mu-1}$ for some $\delta>0$. By definition there exists therefore a dense subset $V\subset U\setminus\{0\}$ such that $\tilde{f}_t=\tilde{F}(\cdot,t)$ has $\mu(\tilde{f}_{\mathbb{C}},0)$ critical points $\tilde{p}_1,\dots,\tilde{p}_{\mu}$ with distinct critical values, for any $t\in V$. Let 
$$F: M\times U\to\mathbb{R},$$
$$F(x_1,\dots,x_{n+1})=\tilde{F}(x_1,x_2)+\sum_{i=3}^{n+1} \epsilon_i x_i^2.$$
Then for each $t\in V$ the critical points of $f_t=F(\cdot, t)$ are
$$p_i=\tilde{p}_i\times (0,\dots,0)\subset M,$$
where $i=1,\dots,\mu$. Therefore $F$ is a strong $\mathbb{R}$-morsification.
\end{item}

\end{enumerate}
\end{proof}

\section{Preparatory Lemmas}
Consider a germ of polynomial maps $f: (\mathbb{R}^{n+1},0)\to(\mathbb{R},0)$ with an isolated singular point in the origin. Let $f: M\to\mathbb{R}$ be a representative of this germ in a neighborhood $M\subset\mathbb{R}^{n+1}$ of the origin.

\begin{theorem}[{see e.g \cite[Theorem 4.2]{Seade}}]\label{Milnor fibrations} There exists $\delta_0>0$ such that for any $\delta\in (0,\delta_0)$ there exists $\epsilon_0=\epsilon_0(\delta)>0$ such that for any $\epsilon\in (0, \epsilon_0)$,
$$f: f^{-1}((0,\epsilon])\cap \bar{\mathbb{B}}_{\delta}\to (0,\epsilon],$$
$$f: f^{-1}((0,-\epsilon])\cap \bar{\mathbb{B}}_{\delta}\to (0,-\epsilon].$$
are the projections of trivial topological fibrations, where $\bar{\mathbb{B}}_{\delta}\subset\mathbb{R}^{n+1}$ denotes the closed ball of radius $\delta$ centered at the origin. Moreover each fibre of these fibrations are smooth manifolds with boundary.
\end{theorem}

Let us fix Milnor data $(\epsilon_0,\delta_0)$ as in the Theorem \hyperref[Milnor fibrations]{\ref*{Milnor fibrations}} and let $\epsilon\in (0,\epsilon_0)$. For any $\eta\in (0, \epsilon]$ let us write 
$$\bar{\mathcal{F}}_{\eta}^{+}=f^{-1}(\eta)\cap\bar{\mathbb{B}}_{\delta},\qquad \bar{\mathcal{F}}_{\eta}^{-}=f^{-1}(-\eta)\cap\bar{\mathbb{B}}_{\delta}$$
for the positive and negative Milnor fibres. It follows from the triviality of the fibrations above that if $\tilde{\eta}$ is another real number such that $\tilde{\eta}\in (0,\epsilon]$ then there are homeomorphisms 
$$\bar{\mathcal{F}}_{\eta}^{+}\cong \bar{\mathcal{F}}_{\tilde{\eta}}^{+},\qquad \bar{\mathcal{F}}_{\eta}^{-}\cong \bar{\mathcal{F}}_{\tilde{\eta}}^{-}$$
Recall the first of the Isotopy lemmas of Réné Thom: 
\begin{lemma}[{``Thom's First Isotopy Lemma'' \cite[Proposition 11.1]{Mather}}]\label{Thom} Let $g: \mathcal{M}\to \mathcal{N}$ be a smooth function of smooth manifolds and let $X\subset \mathcal{M}$ be a closed Whitney stratified subset. If $g_{|X}: X\to \mathcal{N}$ is proper and if for any stratum $S\subset X$ the restriction $g_{|S}: S\to \mathcal{N}$ is a submersion then $f_{|X}: X\to f(X)$ is a locally trivial topological fibration. 
\end{lemma}

We now can prove

\begin{lemma}\label{first lemma} Let $f: (\mathbb{R}^{n+1},0)\to(\mathbb{R},0)$ be a germ of polynomial maps defining an isolated singular point in the origin. Let $f: M\to\mathbb{R}$ be a representative on a neighborhood of the origin and let $(\epsilon_0,\delta_0)$ be Milnor data. Suppose that
$$F: M\times [-1,1]^k\to\mathbb{R},\qquad F(\cdot, t)=f_t$$
is a polynomial map such that $F(x, 0)=f(x)$ for all $x\in M$. Then for any $\delta\in (0,\delta_0)$ and any $\eta\in (0, \epsilon(\delta))$ there exists $t_0'=t_0'(\eta)\in \mathbb{R}$ such that for any $t\in [-1,1]^k$ such that $|t|\leq t_0'$ there exists a homeomorphism 
$$\phi: \bar{\mathcal{F}}_{\eta}^{\pm}\to f_t^{-1}(\pm\eta)\cap\bar{\mathbb{B}}_{\delta}.$$
\end{lemma}
\begin{proof} We will prove the assertion for the positive Milnor fibre, the other case is analogous. By the Curve Selection Lemma \cite[Lemma 3.1]{Milnor} one can choose $\delta_0>0$ such that for any $\delta\in (0,\delta_0)$ the sphere $\mathbb{S}_{\delta}\subset\mathbb{R}^{n+1}$ intersects $f^{-1}(0)$ transversally. Moreover one can choose $\epsilon_0(\delta)>0$ such that any $\eta\in (0, \epsilon_0)$ is a regular value of $f: \mathbb{B}_{\delta}\to\mathbb{R}$. Let us fix $\delta\in (0,\delta_0)$. Then by continuity, $f^{-1}(\eta)\cap \mathbb{S}_{\delta}$ is transverse for all $\eta\in (0,\epsilon_0)$. Let 
$$\tilde{F}: M\times [-1,1]^k\to \mathbb{R}\times [-1,1]^k,\qquad \tilde{F}(x,t)=(f_t(x), t).$$
and denote by
$$\pi: \mathbb{R}^n\times\mathbb{R}^k\to\mathbb{R}^k,\qquad \pi(x,t)=t$$
the standard projection on the parameter $t$. We will show that there exists $t_0'(\eta)$ such that the restricted map
$$\pi: \tilde{F}^{-1}(\{\eta\}\times (-t_0',t_0')^k)\cap(\bar{\mathbb{B}}_{\delta}\times \mathbb{R}^k)\to\mathbb{R}^k$$
is a proper stratified submersion and then use the Isotopy Lemma to obtain the result. To achieve this we will first find a Whitney stratification of the domain.
\begin{enumerate}
\item\label{item one} First note that since any $\eta\in (0,\epsilon_0)$ is a regular value of $f: \mathbb{B}_{\delta}\to\mathbb{R}$ the differential $d\tilde{F}_{| M\times (-1,1)^k}(x,0)$ has maximal rank $k+1$ for any $x\in f^{-1}(\eta)\cap \mathbb{B}_{\delta}$. Therefore, we can assume by continuity that if $\eta$ is fixed then
$$\text{rank }(d\tilde{F}_{|M\times (-1,1)^k})(x,t)=k+1$$
for any $x\in f_t^{-1}(\eta)\cap\mathbb{B}_{\delta}$ and any $t\in (-t_0, t_0)^k$. This means that 
$$\tilde{F}: \mathbb{B}_{\delta}\times (-t_0,t_0)^k\to \mathbb{R}\times \mathbb{R}^k$$
has no critical values in $\{\eta\}\times (-t_0, t_0)^k$. 

\item\label{item two} In particular,
$$\mathscr{T}_1:=\tilde{F}^{-1}(\{\eta\}\times (-t_0,t_0)^k)\cap\left(\mathbb{B}_{\delta}\times (-t_0, t_0)^k\right)$$
is a smooth manifold.
\item\label{item three} We claim that there exists $\tilde{t}_0=\tilde{t}_0(\eta)\in (0,1]$ such that 
$$\mathscr{T}_2:=\tilde{F}^{-1}(\{\eta\}\times (-\tilde{t}_0,\tilde{t}_0)^k)\cap\left(\mathbb{S}_{\delta}\times (-\tilde{t}_0, \tilde{t}_0)^k\right)$$
is a smooth manifold. For this it suffices to show the following claim: the intersection
$$\tilde{F}^{-1}(\{\eta\}\times (-\tilde{t}_0, \tilde{t}_0)^k)\cap \mathbb{S}_{\delta}\times (-\tilde{t}_0, \tilde{t}_0)^k$$
is transverse. This will be shown using the fact that transversality is an open condition. We have that the set of $(x,t)$ with $x\in f_t^{-1}(\eta)$ such that $f_t^{-1}(\eta)\cap \mathbb{S}_{\delta}$ is not transverse at $x$ is precisely the set of $(x,t)$ such that 
$$T_{(x,t)} f_t^{-1}(\eta)\perp v(x)$$
where $v(x)=(x,0)$. Since this is a closed condition, the set where $T_{(x,t)} f_t^{-1}(\eta)$ is not orthogonal to $v(x)$ is an open condition. Since the intersection
$$f^{-1}(\eta)\cap \mathbb{S}_{\delta}$$
is transverse, $T_{(x,0)} f_t^{-1}(\eta)$ is not orthogonal to $v(x)$. Therefore, by openness, there exists a neighborhood of $t=0$ in $\mathbb{R}^k$ such that for any $t$ in that neighborhood and any $x\in f_t^{-1}(\eta)$, the intersection $f_t^{-1}(\eta)\cap\mathbb{S}_{\delta}$ is transverse. Thus there exists a $\tilde{t}_0\in\mathbb{R}$ as claimed. 

\item\label{item four} Let $t_0'=\min(t_0, \tilde{t}_0)$ and consider 
$$\pi_{|\mathscr{T}_i}: \mathscr{T}_i\to (-t_0', t_0')^k,\qquad i=1,2$$
We claim that up to replacing $t_0'$ by a smaller value one has that $\pi_{|\mathscr{T}_i}$ is a submersion, for $i=1,2$. Suppose otherwise and that for any $t_0'$ there exists a point $(x,t)\in\mathscr{T}_1$ in which $\pi_{|\mathscr{T}_1}$ has a critical point. By definition this means that
$$T_{(x,t)} \mathscr{T}_1\perp T_t \mathbb{R}^k\qquad (\ast)$$
and that $d\tilde{F}(x,t)$ is parallel to $T_t\mathbb{R}^k$. Since $(\ast)$ is a closed condition, if $(x_t, t)\in\mathscr{T}_1$ is a sequence such that $|t|\to 0$ then there would exist a limit point $(\tilde{x}, 0)$ with $\tilde{x}\in f^{-1}(-\eta)\cap \bar{\mathbb{B}}_{\delta}$ such that $d\tilde{F}(\tilde{x}, 0)$ is parallel to $T_0\mathbb{R}^k$. But then $d f(\tilde{x})$ would be zero which is impossible since $\tilde{x}\in f^{-1}(\eta)$. Therefore $\pi_{|\mathscr{T}_1}$ is a submersion. In the same vein suppose that $\pi_{|\mathscr{T}_2}$ has a critical point $(x,t)$ which is to say 
$$T_{(x,t)} \mathscr{T}_2\perp T_t \mathbb{R}^k$$
Letting $(x_t, t)\in\mathscr{T}_2$ be a sequence such that $|t|\to 0$ there would exist a limit point $(\tilde{x},0)$ with $\tilde{x}\in f^{-1}(\eta)\cap\mathbb{S}_{\delta}$ such that $d\tilde{F}(\tilde{x},0)$ is parallel to $T_0\mathbb{R}^k$ hence $d f(\tilde{x})=0$ which is impossible since $f^{-1}(\eta)\cap \mathbb{S}_{\delta}$ is transverse at $\tilde{x}$. 

\item Therefore if  
$$\mathscr{W}:=\mathscr{T}_1\cup\mathscr{T}_2$$
then by the the steps \hyperref[item two]{\ref*{item two}} and \hyperref[item three]{\ref*{item three}} above $\mathscr{W}$ is a smooth manifold with boundary hence $\mathcal{T}=\{\mathscr{T}_1, \mathscr{T}_2\}$ is a Whitney stratification. Moreover, by the step \hyperref[item four]{\ref*{item four}}
$$\pi_{|\mathscr{W}}: \mathscr{W}\to (-t_0', t_0')^k$$
is a proper stratified submersion.
\item By Thom's first Isotopy Lemma \hyperref[Thom]{\ref*{Thom}} (see \cite[Proposition 11.1]{Mather}) $\pi_{|\mathscr{W}}$ is the projection of a locally trivial fibration. As the base is contractible it is a trivial fibration. The fibers of a trivial fibration are homeomorphic whence it follows that for any $t\in (-t_0', t_0')^k$ there exists a homeomorphism
$$\phi_t: \pi^{-1}_{|\mathscr{W}}(0)\to \pi^{-1}_{|\mathscr{W}}(t),$$
or in other words,
$$\bar{\mathcal{F}}_{\eta}^{+}\cong f_t^{-1}(\eta)\cap\bar{\mathbb{B}}_{\delta}$$
which is what was to be proven.
\end{enumerate}
\end{proof}

\begin{remark}\label{remark one} Restricted to each strata of the fibers the homeomorphism $\phi_t$ is in fact a diffeomorphism.
\end{remark}

\begin{remark} The set
$$f_t^{-1}(\eta)\cap\bar{\mathbb{B}}_{\delta}$$
is \emph{not} a Milnor fibre of $f_t: \mathbb{R}^{n+1}\to\mathbb{R}$ since the function $f_t$ has in general critical points in $\bar{\mathbb{B}}_{\delta}\setminus\{0\}$. 
\end{remark}

Applying this lemma in the case where $F: M\times U \to\mathbb{R}$ is a $\mathbb{R}$-morsification we get the following corollary:

\begin{corollary}\label{approximation corollary} Let $f: (\mathbb{R}^{n+1},0)\to(\mathbb{R},0)$ be a germ of polynomial maps defining an isolated singular point in the origin and let $F: M\times U\to\mathbb{R}$ be a $\mathbb{R}$-morsification. Then there exists $\delta_0>0$ such that for any $\delta\in (0,\delta_0)$ there exists $\epsilon_0>0$ such that for any $\eta\in (0,\epsilon_0)$ there exists $t_0'=t_0'(\eta)\in \mathbb{R}$ such that for any $t\in U\cap\{|t|\leq t_0'\}$ there exists a homeomorphism 
$$\phi: \bar{\mathcal{F}}_{\eta}^{\pm}\to f_t^{-1}(\pm\eta)\cap\bar{\mathbb{B}}_{\delta}.$$
\end{corollary}
Let us denote
$$\mathcal{X}_{t,\eta}=f_t^{-1}(-\eta,\eta)\cap\bar{\mathbb{B}}_{\delta},$$
$$\bar{\mathcal{F}}_{t,\eta}^{+}=f_t^{-1}(\eta)\cap\bar{\mathbb{B}}_{\delta},$$
$$\bar{\mathcal{F}}_{t,\eta}^{-}=f_t^{-1}(-\eta)\cap\bar{\mathbb{B}}_{\delta}.$$

\begin{remark} We have that $\bar{\mathcal{F}}_{t,\eta}^{+}$ and $\bar{\mathcal{F}}_{t,\eta}^{-}$ belong to the boundary of $\bar{\mathcal{X}}_{t,\eta}$ and that the map 
$$f_t: \bar{\mathcal{X}}_{t,\eta}\to\mathbb{R}$$
has the property that its restriction to $\bar{\mathcal{F}}_{t,\eta}^{+}$ and to $\bar{\mathcal{F}}_{t,\eta}^{-}$ respectively has no critical points, by the proof of Lemma \hyperref[first lemma]{\ref*{first lemma}}. 
\end{remark}

\begin{lemma}\label{manifold lemma} For any $\eta\in (0, \epsilon_0(\delta))$ there exists a real number $t_0''>0$ such that for any $t\in U$ with $|t|\leq t_0''$, the subset $\mathcal{X}_{t,\eta}\subset M$ is a smooth manifold with boundary.
\end{lemma}
\begin{proof} Let $\eta\in (0, \epsilon_0(\delta))$ be fixed. By the Lemma \hyperref[first lemma]{\ref*{first lemma}} there exists $t_0'=t_0'(\eta)$ such that for any $t\in U$ such that $|t|\leq t_0'$ one has that $\eta$ is not a critical value of $f_t: \mathbb{B}_{\delta}\to\mathbb{R}$. Now the open sets
$$f_t^{-1}((-\infty, \eta)\cap\mathbb{B}_{\delta},\qquad f_t^{-1}(-\eta, \infty)\cap\mathbb{B}_{\delta}$$
are smooth manifolds for any $t\in U$ such that $|t|\leq t_0'$. Transversality is an open condition so there exists $t_0''\in\mathbb{R}$ such that the intersection 
$$f_t^{-1}(-\eta,\eta)\cap\mathbb{B}_{\delta}$$
is transverse, for $t\in U\cap\{|t|\leq t_0''\}$ hence is a smooth manifold. The same arguments applied to $f_t: \mathbb{S}_{\delta}\to\mathbb{R}$ shows that 
$$f_t^{-1}(-\eta,\eta)\cap\mathbb{S}_{\delta}$$
is a smooth manifold. This proves the assertion.
\end{proof}

We now show that $\bar{\mathcal{X}}_{t,\eta}$ is contractible.

\begin{lemma}\label{main lemma one} Suppose that 
$$F: M\times U\to\mathbb{R}$$ 
is a $\mathbb{R}$-morsification of $f$ where $U\subset [-1,1]^k$ is contractible. Then for any $\eta\in (0, \epsilon_0(\delta))$ there exists $t_0\in\mathbb{R}$ such that for any $t\in U$ with $|t|\leq t_0$ ,
$$\bar{\mathcal{X}}_{t,\eta}=f_{t}^{-1}([-\eta,\eta])\cap\bar{\mathbb{B}}_{\delta}$$
is contractible.
\end{lemma}
\begin{proof} We put $t_0=\min(t_0'(\eta), t_0''(\eta))$ in the notation of Lemma \hyperref[first lemma]{\ref*{first lemma}} and \hyperref[manifold lemma]{\ref*{manifold lemma}}. Let
$$\tilde{F}: M\times U\to \mathbb{R}\times U,\qquad \tilde{F}(x,t)=(f_t(x),t)$$
and let 
$$\pi: \mathbb{R}^{n+1}\times\mathbb{R}^k\to\mathbb{R}^k,\qquad \pi(x,t)=t$$ 
be the standard projection. Consider
$$\mathscr{T}_1=\tilde{F}^{-1}((-\eta,\eta)\times U\cap\{|t|<t_0\}))\cap \mathbb{B}_{\delta}\times U,$$
$$\mathscr{T}_2=\tilde{F}^{-1}((-\eta,\eta)\times U\cap\{|t|<t_0\}))\cap \mathbb{S}_{\delta}\times U,$$
$$\mathscr{T}_3^{\pm}=\tilde{F}^{-1}(\{\pm \eta\}\times U\cap\{|t|<t_0\}))\cap \mathbb{B}_{\delta}\times U,$$
$$\mathscr{T}_4^{\pm}=\tilde{F}^{-1}(\{\pm \eta\}\times U\cap\{|t|<t_0\}))\cap \mathbb{S}_{\delta}\times U,$$
Then the $\mathscr{T}_1$ and $\mathscr{T}_2$ are smooth manifolds and by the proof of the Lemma \hyperref[first lemma]{\ref*{first lemma}} $\mathscr{T}_3^{\pm}$ and $\mathscr{T}_4^{\pm}$ are smooth manifolds. Therefore if 
$$\mathscr{Y}=\tilde{F}^{-1}([-\eta,\eta]\times U\cap\{|t|< t_0\})\cap(\bar{\mathbb{B}}_{\delta}\times U)$$
then 
$$\mathcal{T}=\{\mathscr{T}_1,\mathscr{T}_2,\mathscr{T}_3^{\pm},\mathscr{T}_4^{\pm}\}$$
is a Whitney stratification of $\mathscr{Y}$ by transversality of their intersections. We will prove that the restriction of the projection
$$\pi: \mathscr{Y}\to U\cap\{|t|<t_0\}$$
is a proper stratified submersion. First of all it is a proper map since any fiber is compact, of the form 
$$\pi_{|\mathscr{Y}}^{-1}(s)=f_t^{-1}([-\eta,\eta])\cap \bar{\mathbb{B}}_{\delta}$$
By the proof of Lemma \hyperref[first lemma]{\ref*{first lemma}}
$$\pi_{|\mathscr{T}_i^{\pm}}: \mathscr{T}_i^{\pm}\to U\cap\{|t|<t_0\},\qquad i=3,4$$
are submersions. Now since 
$$\dim \mathscr{T}_1=\dim (M\times U)=n+k+1$$
is of maximal dimension and since the stratum $\mathscr{T}_1$ is a smooth manifold and in particular open, $\pi_{|\mathscr{T}_1}$ is a submersion and since the projection $\pi_{|\mathbb{B}_{\delta}\times U}$ is a submersion. It thus remains to show that $\pi_{|\mathscr{T}_2}$ is a submersion. Again this follows since 
$$\dim \mathscr{T}_2=\dim (M\times U)-1=n+k,$$
since $\mathscr{T}_2$ is smooth manifold and in particular open and since the projection $\pi_{|\mathbb{S}_{\delta}\times U}$ is a submersion. Therefore one has that 
$$\pi_{|\mathscr{Y}}: \mathscr{Y}\to U\cap\{|t|<t_0\}$$ 
is a stratified proper submersion. By the First Isotopy Lemma of Thom \hyperref[Thom]{\ref*{Thom}} (see \cite[Proposition 11.1]{Mather}) it is therefore the projection of a locally trivial fibration, and as the codomain is contractible by assumption, it is a trivial fibration. In particular its fibres are homeomorphic thus 
$$\pi_{|\mathscr{Y}}^{-1}(t)=\bar{\mathcal{X}}_{t,\eta}\cong f^{-1}([-\eta,\eta])\cap\bar{\mathbb{B}}_{\delta}=\pi_{|\mathscr{Y}}^{-1}(0)$$

Therefore, to finish the proof it suffices to show that $f^{-1}([-\eta,\eta])\cap\bar{\mathbb{B}}_{\delta}$ is contractible. But it follows from (\cite[Proposition 1.6]{durfee}) that the inclusion $f^{-1}(0)\cap\bar{\mathbb{B}}_{\delta}\subset f^{-1}([-\eta,\eta])\cap\bar{\mathbb{B}}_{\delta}$ is a homotopy equivalence. By the Local Conic Structure of algebraic sets $f^{-1}(0)\cap\bar{\mathbb{B}}_{\delta}$ is a cone over its boundary hence contractible. Hence $f^{-1}([-\eta,\eta])\cap\bar{\mathbb{B}}_{\delta}$ and therefore also $\bar{\mathcal{X}}_{t,\eta}$ is contractible.
\end{proof}

\newpage
\section{On the Homotopy Type of The Real Milnor Fibres}\label{section 49}
\subsection{Introduction}
The idea now is to Whitney stratify $\bar{\mathcal{X}}_{t,\eta}$ and to use the nonproper version of Stratified Morse Theory to obtain the homotopy type of the real Milnor fibres. More precisely we will in Theorem \hyperref[main theorem one]{\ref*{main theorem one}} show that a tubular neighborhood of a real Milnor fibre becomes homeomorphic to a contractible space after succesively adjoining handles attached via embeddings, with each handle corresponding to a critical point of an $\mathbb{R}$-morsification.\\
\subsection{The Situation}
Let $f: (\mathbb{R}^{n+1}, 0)\to (\mathbb{R}, 0)$ be a germ of isolated singularity and let 
$$F: M\times U\to\mathbb{R},\qquad f_t=F(\cdot, t)$$
be a $\mathbb{R}$-morsification of this germ, where $M\subset\mathbb{R}^{n+1}$ is a neighborhood of the origin. Fix Milnor data $(\epsilon_0, \delta_0)$ at the origin for $f: M\to \mathbb{R}$, as in Theorem \hyperref[Milnor fibrations]{\ref*{Milnor fibrations}}.\\

We shall only consider those critical points which are bounded in the following sense.

\begin{definition}\label{number} Let $\delta\in(0,\delta_0)$ and $\epsilon\in(0,\delta)$ be fixed. For any $t\in V$ with $V\subset U\setminus\{0\}$ as in Definition \hyperref[definition morsification]{\ref*{definition morsification}} let $p_1,\dots, p_{m(t)}$ denote the critical points of $f_t: M\to\mathbb{R}$ such that $p_i\in\bar{\mathbb{B}}_{\delta}$ for all $i=1,\dots, m(t)$ and let $s_i=f_t(p_i)$ denote their values.
\end{definition}

\begin{remark}\label{bounded remark} We claim that there exists $t_0\in\mathbb{R}$ such that $s_i(t)\in (-\epsilon,\epsilon)$ for all $t\in V\cap\{|t|\leq t_0\}$. Write $U'=U\cap\{|t|\leq t_0\}$. We first claim that if $\mathscr{C}_t$ denotes the set of $t\in U'$ such that $p_i(t)\in \bar{\mathbb{B}}_{\delta}$ is a critical point of $f_t: M\to \mathbb{R}$ then $\mathscr{C}_t$ is closed in $U$. Consider 
$$\mathscr{C}=\{(x,t)\in \bar{\mathbb{B}}_{\delta}\times U'\ |\ \text{rank }\text{Jac}(f_t)(x)=0\}$$
Then $\mathscr{C}$ is closed in the compact set $\bar{\mathbb{B}}_{\delta}\times U'\subset M\times \mathbb{R}^k$ with respect to the Euclidean topology and so if 
$$\pi: \mathbb{R}^{n+1}\times \mathbb{R}^k\to \mathbb{R}^k,\qquad \pi(x,t)=t$$
denotes the standard projection onto the second factor then the image $\mathscr{C}_t=\pi(\mathscr{C})$ is closed in $U'$. Therefore, by continuity of $t\mapsto df_t$, if 
$$(p_{t_n})_{n\in\mathbb{N}},\qquad |p_{t_n}|\leq \delta$$
is a bounded sequence of critical points such that $t_n\in U'$ and $\lim_{n\to \infty}|t_n|=0$ then the limit $p=p_0$ exists, is bounded $|p_0|\leq \delta$, and is a critical point of $f$. As the origin is the unique critical point of $f$ in the ball $\bar{\mathbb{B}}_{\delta}\subset M$ it follows that 
$$\forall p_i(t)\in \bar{\mathbb{B}}_{\delta},\qquad \lim_{|t|\to 0} p_i(t)=0$$
Since $f_t: M\to \mathbb{R}$ is a representative of a germ $f_t: (\mathbb{R}^{n+1}, 0)\to (\mathbb{R}, 0)$ one has $f_t(0)=0$ hence the critical values satisfy $s_i(t)\in (-\eta,\eta)$ for all $i=1,\dots, m(t)$ and for all $t\in V\cap\{|t|\leq t_0\}$ for $t_0$ sufficiently small. 
\end{remark}

\subsection{The Theorem}
The following theorem concerns the homotopy type of the real Milnor fibres of $f$ at the origin. In its statement we shall for notational simplicity write $\mathcal{N}(\bar{\mathcal{F}}_{\eta}^{\pm})$ for the cartesian product of $(0,1)$ with the positive, respectively negative, Milnor fibres.

\begin{theorem}\label{main theorem one} Suppose that $F: M\times U\to\mathbb{R}$ is a $\mathbb{R}$-morsification. For any $\eta\in (0, \epsilon(\delta))$ there exists $t_0=t_0(\eta)\in \mathbb{R}$ such that the following holds. Let $t\in V$ (where $V$ is as in Definition \hyperref[definition morsification]{\ref*{definition morsification}}) such that $|t|\leq t_0$ be fixed and let $p_1,\dots,p_m\in\bar{\mathbb{B}}_{\delta}$ denote the critical points of $f_t: M\to \mathbb{R}$ lying inside the ball of radius $\delta$ centered at the origin and let $\lambda(p_i)$ denote their indices. Then $\mathcal{X}_{t,\eta}$ is contractible and there exist embeddings
\[h_i^{+}: \{
    \begin{array}{ll}
      \partial(\mathbb{D}^{\lambda(p_i)}\times \partial \mathbb{D}^{n+1-\lambda(p_i)})\to \mathcal{N}(\bar{\mathcal{F}}_{\eta}^{+})\cup \bigcup_{\substack{1\leq j\leq i-1 \\ h_j^{+}}} \left(\mathbb{D}^{\lambda(p_j)}\times \mathbb{D}^{n+1-\lambda(p_j)}\right),\quad i\geq 2\\
      \partial(\mathbb{D}^{\lambda(p_1)}\times \partial \mathbb{D}^{n+1-\lambda(p_1)})\to \mathcal{N}(\bar{\mathcal{F}}_{\eta}^{+}),\quad i=1.
\end{array}
\]
\[h_i^{-}: \{
    \begin{array}{ll}
      \partial(\mathbb{D}^{n+1-\lambda(p_i)}\times \partial \mathbb{D}^{\lambda(p_i)})\to \mathcal{N}(\bar{\mathcal{F}}_{\eta}^{-})\cup \bigcup_{\substack{1\leq j\leq i-1 \\ h_j^{-}}} \left(\mathbb{D}^{n+1-\lambda(p_j)}\times \mathbb{D}^{\lambda(p_j)}\right),\quad i\geq 2\\
      \partial(\mathbb{D}^{n+1-\lambda(p_1)}\times \partial \mathbb{D}^{\lambda(p_1)})\to \mathcal{N}(\bar{\mathcal{F}}_{\eta}^{-}),\quad i=1.
\end{array}
\]
and homeomorphism
$$\mathcal{X}_{t,\eta}\cong \mathcal{N}(\bar{\mathcal{F}}_{\eta}^{+})\cup \bigcup_{\substack{1\leq i\leq m \\ h_i^{+}}} \left(\mathbb{D}^{\lambda(p_i)}\times \mathbb{D}^{n+1-\lambda(p_i)}\right)$$
$$\mathcal{X}_{t,\eta}\cong \mathcal{N}(\bar{\mathcal{F}}_{\eta}^{-})\cup \bigcup_{\substack{1\leq i\leq m \\ h_i^{-}}} \left(\mathbb{D}^{n+1-\lambda(p_i)}\times \mathbb{D}^{\lambda(p_i)}\right).$$
where each handle $\mathbb{D}^{\lambda(p_i)}\times \mathbb{D}^{n+1-\lambda(p_i)}$ (respectively $\mathbb{D}^{n+1-\lambda(p_i)}\times \mathbb{D}^{\lambda(p_i)}$) is attached along $\mathbb{D}^{\lambda(p_i)}\times \partial \mathbb{D}^{n+1-\lambda(p_i)}$ (respectively along $\mathbb{D}^{n+1-\lambda(p_i)}\times \partial \mathbb{D}^{\lambda(p_i)}$) via $h_i^{+}$ (respectively via $h_i^{-}$).
\end{theorem}
\begin{proof} Let $t_0(\eta)=\min(t_0'(\eta), t_0''(\eta))$ in the notation of Lemma \hyperref[main lemma one]{\ref*{main lemma one}}. Then for any $t\in U\cap\{|t|\leq t_0\}$ we have that $\pm\eta$ are regular values of $f_t: \bar{\mathcal{X}}_{t,\eta}\to\mathbb{R}$, and consequently no point of $\bar{\mathcal{F}}_{t,\eta}^{\pm}$ are critical points. We now fix $t\in V$ such that $|t|<t_0$.
\begin{enumerate}
\item Let 
$$\mathscr{S}_1=\text{int} \mathcal{X}_{t,\eta},\quad \mathscr{S}_2=\mathcal{X}_{t,\eta}\cap\mathbb{S}_{\delta}$$
$$\mathscr{S}_3=\mathcal{F}_{t,\eta}^{+},\quad \mathscr{S}_4=\mathcal{F}_{t,\eta}^{-},\quad \mathscr{S}_5=\partial\bar{\mathcal{F}}^{+}_{t,\eta},\quad \mathscr{S}_6=\partial\bar{\mathcal{F}}_{t,\eta}^{-}$$
and let $\mathcal{S}=\{\mathscr{S}_1,\dots\mathscr{S}_6\}$. As $|t|<t_0''$ it follows from the Lemma \hyperref[manifold lemma]{\ref*{manifold lemma}} that $\mathscr{S}_i$ are smooth manifolds, for $i=1,2$. As $|t|<t_0'$ it follows from the steps \hyperref[item one]{\ref*{item one}} and \hyperref[item three]{\ref*{item three}} of the proof of \hyperref[first lemma]{\ref*{first lemma}} that $\mathscr{S}_i$ are smooth manifolds, for $i=3,\dots, 6$. Moreover as 
$$f_t^{-1}([-\eta,\eta])\cap\mathbb{B}_{\delta}=\mathscr{S}_1\cup\mathscr{S}_3\cup\mathscr{S}_4$$
is a manifold with boundary its strata satisfy the Whitney conditions (see e.g \cite[Definition 1.2]{Goresky}). As the same holds for 
$$f_t^{-1}([-\eta,\eta])\cap\mathbb{S}_{\delta}=\mathscr{S}_2\cup \mathscr{S}_5\cup\mathscr{S}_6$$ 
its strata satisfy the Whitney conditions. But then by tranversality (step \hyperref[item three]{\ref*{item three}} in the proof of \hyperref[first lemma]{\ref*{first lemma}}) the strata of $\mathcal{S}$ satisfy also the Whitney conditions. Therefore it is a Whitney stratification of $\bar{\mathcal{X}}_{t,\eta}$ such that $\mathcal{X}_{t,\eta}$ is a union of strata.

\item Let us consider the function $f_t: M\to \mathbb{R}$. It restricts to a proper function
$$f_t: \bar{\mathcal{X}}_{t,\eta}\to [-\eta,\eta]\subset\mathbb{R}.$$
At the points $p\in \bar{\mathcal{F}}_{t,\eta}^{\pm}$ this function have depraved stratified critical points since $f_t$ is constant there. However any critical point $$p\in\bar{\mathcal{X}}_{t,\eta}\setminus\bar{\mathcal{F}}_{t,\eta}^{+}\cup \bar{\mathcal{F}}_{t,\eta}^{-}$$
is necessarily nondepraved, as we now prove. So we need to show that if $p\in\mathscr{S}_i$ is a critical point of 
$$f_{t| \mathscr{S}_i}: \mathscr{S}_i\to \mathbb{R},\qquad i=1,2$$   
then it is nondepraved. For this, we first show $f_{t|\mathscr{S}_2}$ has no critical point. Indeed suppose not. Then 
$$f_{t|\mathbb{S}_{\delta}}: \mathbb{S}_{\delta}\to \mathbb{R}$$
would have a critical point $p\in\mathbb{S}_{\delta}$ with critical value $s\in (-\eta,\eta)$ hence $df_{t|\mathbb{S}_{\delta}}(p)=0$ and 
$$T_p \mathbb{S}_{\delta}\subset \ker df_{t|\mathbb{S}_{\delta}}(p).\qquad(\ast)$$ 
However the intersections $f_t^{-1}(\pm\eta)\cap \mathbb{S}_{\delta}$ are transverse, by Lemma \hyperref[first lemma]{\ref*{first lemma}} for all $t$ sufficiently small (that is, up to replacing $t_0$ be a smaller value). Therefore by continuity if $t$ is fixed, then $f_t^{-1}(s)\cap\mathbb{B}_{\delta}$ is transverse for all $s\in (-\eta,\eta)$. Then $(\ast)$ implies that $\ker df_{t|\mathbb{S}_{\delta}}(p)=T_p f_t^{-1}(s)$ and we conclude that $T_p \mathbb{S}_{\delta}\subset T_p f_t^{-1}(s)$, which is impossible since we had a transverse intersection. Therefore $\mathscr{S}_2$ has no critical points and it remains to show that any critical point $p$ in $\mathscr{S}_1$ is nondepraved. Since $\mathscr{S}_1$ is a smooth manifold one can apply \cite[Proposition 2.4]{Goresky} to 
$$f_{t| \mathscr{S}_1}: \mathscr{S}_1\to\mathbb{R}$$
which gives that any critical point is nondepraved. All said, any stratified critical point of $f_t: \bar{\mathcal{X}}_{t,\eta}\to\mathbb{R}$ which belongs to a stratum $\mathscr{S}_i$ for $i=1,2$ necessarily belong to the first, and is a nondepraved point.

\item\label{last first} We now apply the nonproper version of Stratified Morse Theory with the stratified subset $\mathcal{X}_{t,\eta}$ and the proper function $f_t: \bar{\mathcal{X}}_{t,\eta}\to\mathbb{R}$. We will use Theorem \cite[Theorem 10.4]{Goresky} in conjunction with the Fundamental Theorem \cite[Theorem 10.5]{Goresky}.\\

Let $p\in\mathcal{X}_{t,\eta}$ be a critical point with value $s=f_t(p)$. By the previous step, $p\in\mathscr{S}_1$ is an interior point. As $\mathscr{S}_1$ has maximal dimension $\dim \mathscr{S}_1=\dim M$ it follows that its normal slice is $N_{\mathcal{X}_{t,\eta}}(p)=p$ whence it follows by its very definition that the normal Morse data is 
$$(\mathscr{N}_{\mathcal{X}_{t,\eta}}, \mathscr{N}_{\mathcal{X}_{t,\eta}})=(p, \emptyset)$$
Therefore by the \cite[Theorem 10.5]{Goresky} the local Morse data for $f_{t|\mathcal{X}_{t,\eta}}$ at the point $p\in\mathscr{S}_1$ is homeomorphic to the tangential Morse data 
$$(\mathscr{A}_{\mathcal{X}_{t,\eta}}, \mathscr{B}_{\mathcal{X}_{t,\eta}})\cong (\mathscr{T}_{\mathcal{X}_{t,\eta}}, \mathscr{T}'_{\mathcal{X}_{t,\eta}})$$
for $f_{t|\mathcal{X}_{t,\eta}}$ at the point, which is by definition the tangential Morse data of $f_t$ at $p$. By \cite[Proposition 4.5]{Goresky} the tangential Morse data at $p$ is homeomorphic by a decomposition-preserving homeomorphism to
$$(\mathscr{T}_{\mathcal{X}_{t,\eta}}, \mathscr{T}'_{\mathcal{X}_{t,\eta}})=(\mathbb{D}^{n+1-\lambda_S(p)}\times \mathbb{D}^{\lambda_S(p)}, \mathbb{D}^{n+1-\lambda_S(p)}\times \partial \mathbb{D}^{\lambda_S(p)})$$
Thus by \cite[Theorem 10.4]{Goresky} one has that if $[s-\gamma, s+\gamma]$ contains no critical values except $s=f_t(p)$ then 
$$f_t^{-1}(-\eta, s+\gamma)\cap \bar{\mathbb{B}}_{\delta}\cong f^{-1}(-\eta, s-\gamma)\cap \bar{\mathbb{B}}_{\delta}\cup_{\mathscr{T}'_{\mathcal{X}_{t,\eta}}} \mathscr{T}_{\mathcal{X}_{t,\eta}}$$
which is to say that there exists an embedding 
$$h_s: \mathbb{D}^{n+1-\lambda_S(p)}\times \partial \mathbb{D}^{\lambda_S(p)}\to f_t^{-1}(-\eta, s-\gamma)\cap \bar{\mathbb{B}}_{\delta}$$
such that the identity map extends to a homeomorphism
$$f_t^{-1}(-\eta, s+\gamma)\cap \bar{\mathbb{B}}_{\delta}\cong$$
$$f_t^{-1}(-\eta, s-\gamma)\cap \bar{\mathbb{B}}_{\delta}\cup_h \left(\mathbb{D}^{n+1-\lambda_S(p)}\times \mathbb{D}^{\lambda_S(p)}\right)$$

\item\label{helvete} By the Remark \hyperref[bounded remark]{\ref*{bounded remark}} if $p_1,\dots, p_m\in \mathbb{B}_{\delta}$ denote the critical values of $f_{t|\mathscr{S}_1}$ we can order the corresponding critical values $s_i=f_t(p_i)$ as follows 
$$-\eta<s_1<\dots<s_m<\eta.$$
For each $i=1,\dots, m$ choose a real number $\gamma_i>0$ such that $s_i$ is the only critical value contained in $[s_i-\gamma_i, s_i+\gamma_i]$ and such that
$$s_i+\gamma_i=s_{i+1}-\gamma_{i+1},\qquad i=1,\dots, m-1$$
The previous step \hyperref[last first]{\ref*{last first}} gives that there exist an embedding 
$$h_m^{-}: \mathbb{D}^{n+1-\lambda_S(p_m)}\times \partial \mathbb{D}^{\lambda_S(p_m)}\to f_t^{-1}(-\eta, s_m-\gamma_m)\cap\bar{\mathbb{B}}_{\delta}$$
such that 
$$\mathcal{X}_{t,\eta}\cong f_t^{-1}(-\eta, s_m-\gamma_m)\cap\bar{\mathbb{B}}_{\delta}\cup_{h_m^{-}} \left(\mathbb{D}^{n+1-\lambda_{\mathscr{S}_1}(p_m)}\times \mathbb{D}^{\lambda_{\mathscr{S}_1}(p_m)}\right)$$
Similarly there exists an embedding 
$$h_{m-1}^{-}: \mathbb{D}^{n+1-\lambda_S(p_{m-1})}\times \partial \mathbb{D}^{\lambda_S(p_{m-1})}\to f_t^{-1}(-\eta, s_{m-1}-\gamma_{m-1})\cap\bar{\mathbb{B}}_{\delta}$$
such that 
$$f_t^{-1}(-\eta, s_m-\gamma_m)\cap\bar{\mathbb{B}}_{\delta}=f_t^{-1}(-\eta, s_{m-1}+\gamma_{m-1})\cap\bar{\mathbb{B}}_{\delta}$$
$$\cong f_t^{-1}(-\eta, s_{m-1}-\gamma_{m-1})\cap\bar{\mathbb{B}}_{\delta}\cup_{h_{m-1}^{-}} \left(\mathbb{D}^{n+1-\lambda_{\mathscr{S}_1}(p_{m-1})}\times \mathbb{D}^{\lambda_{\mathscr{S}_1}(p_{m-1})}\right).$$
Therefore $\mathcal{X}_{t,\eta}$ is homeomorphic to
$$\left(f_t^{-1}(-\eta, s_{m-1}-\gamma_{m-1})\cap\bar{\mathbb{B}}_{\delta}\cup_{h_{m-1}^{-}} \left(\mathbb{D}^{n+1-\lambda_{\mathscr{S}_1}(p_{m-1})}\times \mathbb{D}^{\lambda_{\mathscr{S}_1}(p_{m-1})}\right)\right)$$
$$\cup_{h_m^{-}} \left(\mathbb{D}^{n+1-\lambda_{\mathscr{S}_1}(p_m)}\times \mathbb{D}^{\lambda_{\mathscr{S}_1}(p_m)}\right)$$
Continuing inductively this yields 
$$\mathcal{X}_{t,\eta}\cong f_t^{-1}(-\eta, s_1-\gamma_1)\cap\bar{\mathbb{B}}_{\delta}\cup \bigcup_{\substack{1\leq i\leq m \\ h_i^{-}}} \left(\mathbb{D}^{n+1-\lambda_{\mathscr{S}_1}(p_i)}\times \mathbb{D}^{\lambda_{\mathscr{S}_1}(p_i)}\right).$$
\item\label{dra helvete} Since the dimension of the stratum $\mathscr{S}_1$ is maximal, the index of $f_{t|\mathscr{S}_1}$ at a critical point is equal to the index of $f_{t}: M\to\mathbb{R}$ so 
$$\lambda_{\mathscr{S}_1}(p_i)=\lambda(p_i).$$
It therefore remains to show that 
$$f_t^{-1}(-\eta,s_1-\gamma_1)\cap\bar{\mathbb{B}}_{\delta}\cong \bar{\mathcal{F}}_{\eta}^{-}\times (0,1)$$
Now, 
$$f_t^{-1}(-\eta, s_1-\gamma_1)\cap\bar{\mathbb{B}}_{\delta}\cong f_t^{-1}(-\eta+\gamma_1/2)\cap\bar{\mathbb{B}}_{\delta}\times (-\eta, s_1-\gamma_1)$$ 
and 
$$f_t^{-1}(-\eta-\gamma_1,s_1-\gamma_1)\cap\bar{\mathbb{B}}_{\delta}\cong f_t^{-1}(-\eta)\cap\bar{\mathbb{B}}_{\delta}\times (-\eta-\gamma_1,s_1-\gamma_1)$$
by Thom's Isotopy Lemma \hyperref[Thom]{\ref*{Thom}} (\cite[Proposition 11.1]{Mather}) because by construction $f_t: M\to\mathbb{R}$ has no critical values in $(-\eta-\gamma_1 ,s_1-\gamma_1)$ and it moreover follows that 
$$f_t^{-1}(-\eta, s_1-\gamma_1)\cap\bar{\mathbb{B}}_{\delta}\cong f_t^{-1}(-\eta-\gamma_1,s_1-\gamma_1)\cap\bar{\mathbb{B}}_{\delta}.$$
But then 
$$\bar{\mathcal{F}}_{t,\eta}^{-}=f_t^{-1}(-\eta)\cap\bar{\mathbb{B}}_{\delta}\cong f_t^{-1}(-\eta+\gamma_1/2)\cap\bar{\mathbb{B}}_{\delta}$$
and consequently $f_t^{-1}(-\eta,s_1-\gamma_1)\cap\bar{\mathbb{B}}_{\delta}\cong \bar{\mathcal{F}}_{\eta}^{-}\times (0,1)$ because $\bar{\mathcal{F}}_{t,\eta}^{-}\cong \bar{\mathcal{F}}_{\eta}$ by 
Lemma \hyperref[first lemma]{\ref*{first lemma}}. Finally that $\mathcal{X}_{t,\eta}$ is contractible follows from Lemma \hyperref[main lemma one]{\ref*{main lemma one}}.
\end{enumerate}
\end{proof}

\newpage

Figures (\hyperref[figure contractible]{\ref*{figure contractible}}--\hyperref[fag]{\ref*{fag}}) illustrate Theorem \hyperref[main theorem one]{\ref*{main theorem one}} in the case of a quadratic singularity.\\

\begin{figure}[h]

 \centering
  \begin{tikzpicture}[scale=0.4]

\clip (0,0) circle [radius=4.5 cm];

\begin{scope}
\draw[black, opacity=1] (0,0) circle (4.5);

\path[fill=gray,opacity=1] (-10, 10) circle (13.5);
\path[fill=gray,opacity=1] (10, -10) circle (13.5);
\fill[blue, even odd rule, opacity=1] (5,5) circle (4) (5,5) circle (3.9);
\fill[blue, even odd rule, opacity=1] (-5,-5) circle (4) (-5,-5) circle (3.9);

\fill[black, even odd rule, opacity=1] (5,5) circle (2) (5,5) circle (3.9);

\fill[black, even odd rule, opacity=1] (-5,-5) circle (2) (-5,-5) circle (3.9);
\end{scope}


\draw (0,-4.5) -- (0, 4.5);
\draw (-4.5,0) -- (4.5, 0);

\filldraw [black] (0,0) circle (2pt);
\filldraw [black] (0,-2) circle (2pt);
\filldraw [black] (0,-4) circle (2pt);
\filldraw [black] (0, 2) circle (2pt);
\filldraw [black] (0, 4) circle (2pt);
\filldraw [black] (2,0) circle (2pt);
\filldraw [black] (-2,0) circle (2pt);
\filldraw [black] (4,0) circle (2pt);
\filldraw [black] (-4,0) circle (2pt);
\end{tikzpicture}
\caption{The positive Milnor fibre, in blue, of $f=xy$ with the handle $\mathbb{D}^1\times \mathbb{D}^1$, represented by the white strip, attached along the two disjoint segments $\mathbb{D}^1\times\mathbb{S}^0$.}\label{figure contractible}
  \centering
  \begin{tikzpicture}[scale=0.4]

\clip (0,0) circle [radius=4.5 cm];

\begin{scope}

\draw[black, opacity=1] (0,0) circle (4.5);

\fill[blue, opacity=1] (5,5) circle (4) (5,5) circle (3.9);
\fill[blue, opacity=1] (-5,-5) circle (4) (-5,-5) circle (3.9);

\fill[gray, opacity=1] (5,5) circle (2) (5,5) circle (3.9);

\fill[gray, opacity=1] (-5,-5) circle (2) (-5,-5) circle (3.9);

\fill[blue, opacity=1] (5,-5) circle (4) (5,-5) circle (3.9);
\fill[blue, opacity=1] (-5,5) circle (4) (-5,5) circle (3.9);

\fill[gray, opacity=1] (-5,5) circle (2) (-5,5) circle (3.9);

\fill[gray, opacity=1] (5,-5) circle (2) (5,-5) circle (3.9);

\end{scope}


\draw (0,-4.5) -- (0, 4.5);
\draw (-4.5,0) -- (4.5, 0);

\filldraw [black] (0,0) circle (2pt);
\filldraw [black] (0,-2) circle (2pt);
\filldraw [black] (0,-4) circle (2pt);
\filldraw [black] (0, 2) circle (2pt);
\filldraw [black] (0, 4) circle (2pt);
\filldraw [black] (2,0) circle (2pt);
\filldraw [black] (-2,0) circle (2pt);
\filldraw [black] (4,0) circle (2pt);
\filldraw [black] (-4,0) circle (2pt);
\end{tikzpicture}
\caption{The region in white, representing $\mathcal{X}_{t,\eta}$ is homeomorphic to the space obtained from the positive Milnor with a handle attached as in Figure \hyperref[figure contractible]{\ref*{figure contractible}}.}\label{fag}
  
\end{figure}

\begin{corollary}\label{corollary main theorem one} With the same assumptions as in Theorem \hyperref[main theorem one]{\ref*{main theorem one}} there exist attaching maps 
\[H_i^{+}: \left\{
    \begin{array}{ll}
      \partial \mathbb{D}^{n+1-\lambda(p_i)}\to \bar{\mathcal{F}}_{\eta}^{+}\cup \bigcup_{\substack{1\leq j\leq i-1 \\ H_j^{+}}} \mathbb{D}^{n+1-\lambda(p_j)},\quad i\geq 2\\
      \partial \mathbb{D}^{n+1-\lambda(p_1)}\to \bar{\mathcal{F}}_{\eta}^{+},\quad i=1.
\end{array} \right. \]
\[H_i^{-}: \left\{
    \begin{array}{ll}
      \partial \mathbb{D}^{\lambda(p_i)}\to \bar{\mathcal{F}}_{\eta}^{-}\cup \bigcup_{\substack{1\leq j\leq i-1 \\ H_j^{-}}} \mathbb{D}^{\lambda(p_j)},\quad i\geq 2\\
      \partial \mathbb{D}^{\lambda(p_1)}\to \bar{\mathcal{F}}_{\eta}^{-},\quad i=1.
\end{array} \right. \]

and homotopy equivalences 
$$\mathcal{X}_{t,\eta}\sim \bar{\mathcal{F}}^{+}\cup_{H_1^{+}} \mathbb{D}^{n+1-\lambda(p_1)}\cup\dots\cup_{H_m^{+}} \mathbb{D}^{n+1-\lambda(p_m)}$$ 
$$\mathcal{X}_{t,\eta}\sim \bar{\mathcal{F}}^{+}\cup_{H_1^{-}} \mathbb{D}^{\lambda(p_1)}\cup\dots\cup_{H_m^{-}} \mathbb{D}^{\lambda(p_m)}.$$ 
where each disc $\mathbb{D}^{n+1-\lambda(p_i)}$ (respectively $\mathbb{D}^{\lambda(p_i)}$) is attached along its boundary via $H_i^{+}$ (respectively via $H_i^{-}$).
\end{corollary}
\begin{proof} The proof consists of repeating steps (\hyperref[last first]{\ref*{last first}}--\hyperref[dra helvete]{\ref*{dra helvete}}) of the previous Theorem \hyperref[main theorem one]{\ref*{main theorem one}} but using homotopy Morse data \cite[Definition 3.3]{Goresky} instead of Morse data.\\

Let $t\in V\cap\{|t|\leq t_0\}$ be fixed, with $t_0$ is as Lemma \hyperref[main lemma one]{\ref*{main lemma one}}. In particular $\bar{\mathcal{F}}_{t,\eta}^{-}\cong \bar{\mathcal{F}}_{\eta}$ and $\mathcal{X}_{t,\eta}$ is contractible. Let again $p\in\mathcal{X}_{t,\eta}$ be a critical point of $f_t: \bar{\mathcal{X}}_{t,\eta}\to \mathbb{R}$ and let $s=f_t(p)$. In step \hyperref[last first]{\ref*{last first}} of the proof of Theorem \hyperref[main theorem one]{\ref*{main theorem one}} was established that the local Morse data of $f_{t|\mathcal{X}_{t,\eta}}$ at $p$ is homeomorphic to the tangential Morse data of $f_t$ at $p$. Furthermore it was established that the tangential Morse data of $f_t$ at $p$ is homeomorphic to
$$(\mathbb{D}^{n+1-\lambda_S (p)}\times \mathbb{D}^{\lambda_S (p)}, \mathbb{D}^{n+1-\lambda_S (p)}\times\partial\mathbb{D}^{\lambda_S (p)}).$$
By \cite[Remark 3.5.4]{Goresky} there exists a homotopy equivalence
$$(\mathbb{D}^{n+1-\lambda_S (p)}\times \mathbb{D}^{\lambda_S (p)}, \mathbb{D}^{n+1-\lambda_S (p)}\times\partial\mathbb{D}^{\lambda_S (p)})\sim (\mathbb{D}^{\lambda_S(p)}, \partial \mathbb{D}^{\lambda_S(p)})$$
By \cite[Remark 3.3]{Goresky} this implies that if $\gamma=\gamma(p)>0$ is chosen as in step \hyperref[last first]{\ref*{last first}} of the proof of Theorem \hyperref[main theorem one]{\ref*{main theorem one}} then there exists an attaching map 
$$H_s: \partial\mathbb{D}^{\lambda(p)}\to f_t^{-1}(-\eta, s-\gamma)\cap \bar{\mathbb{B}}_{\delta}$$
and a homotopy equivalence 
$$f_t^{-1}(-\eta, s+\gamma)\cap\bar{\mathbb{B}}_{\delta}\sim f_t^{-1}(-\eta, s-\gamma)\cap\bar{\mathbb{B}}_{\delta}\cup_{H_s} \mathbb{D}^{\lambda(p)}.$$
\begin{enumerate}
 \item In general if real numbers $\gamma_1,\dots, \gamma_m>0$ are chosen as in step \hyperref[helvete]{\ref*{helvete}} of Theorem \hyperref[main theorem one]{\ref*{main theorem one}} then it follows that there exist for each $i=1,\dots, m$ attaching maps
$$\tilde{H}_i^{-}: \partial \mathbb{D}^{\lambda(p_i)}\to f_t^{-1}(-\eta, s_i-\gamma_i)\cap\bar{\mathbb{B}}_{\delta}$$
and homotopy equivalences 
$$f_t^{-1}(-\eta, s_i+\gamma_i)\cap\bar{\mathbb{B}}_{\delta}\sim f_t^{-1}(-\eta, s_i-\gamma_i)\cap\bar{\mathbb{B}}_{\delta}\cup_{\tilde{H}_i^{-}} \mathbb{D}^{\lambda(p_i)}.$$
This gives 
$$\mathcal{X}_{t,\eta}\sim f_t^{-1}(-\eta, s_m-\gamma_m)\cap \bar{\mathbb{B}}_{\delta}\cup_{\tilde{H}_m^{-}} \mathbb{D}^{\lambda(p_m)}$$
$$\sim f_t^{-1}(-\eta, s_{m-1}-\gamma_{m-1})\cap \bar{\mathbb{B}}_{\delta}\cup_{\tilde{H}_{m-1}^{-}} \mathbb{D}^{\lambda(p_{m-1})}\cup_{\tilde{H}_m^{-}} \mathbb{D}^{\lambda(p_m)}.$$
Continuing in this manner one obtains 
$$\mathcal{X}_{t,\eta}\sim f_t^{-1}(-\eta, s_1-\gamma_1)\cap\bar{\mathbb{B}}_{\delta}\cup_{\tilde{H}_1^{-}} \mathbb{D}^{\lambda(p_1)}\cup\dots \cup_{\tilde{H}_m^{-}} \mathbb{D}^{\lambda(p_m)}.$$
Since by step \hyperref[dra helvete]{\ref*{dra helvete}} of the proof of Theorem \hyperref[main theorem one]{\ref*{main theorem one}}, 
$$f_t^{-1}(-\eta, s_1-\gamma_1)\cap\bar{\mathbb{B}}_{\delta}\cong \bar{\mathcal{F}}_{\eta}^{-}\times (0,1)$$
one gets that
$$\mathcal{X}_{t,\eta}\sim \bar{\mathcal{F}}^{-}\times (0,1)\cup_{\tilde{H}_1^{-}} \mathbb{D}^{\lambda(p_1)}\cup\dots\cup_{\tilde{H}_m^{-}} \mathbb{D}^{\lambda(p_m)}.$$
\item Let
$$\phi_1: \bar{\mathcal{F}}^{-}\times (0,1)\to \bar{\mathcal{F}}^{-}$$
be a deformation retraction. Then $\phi_1$ extends by \cite[Lemma 3.7]{Morse} to a homotopy equivalence 
$$\phi_2: \bar{\mathcal{F}}^{-}\times (0,1)\cup_{\tilde{H}_1^{-}} \mathbb{D}^{\lambda(p_1)}\to \bar{\mathcal{F}}^{-}\cup_{\phi_1\circ \tilde{H}_1^{-}} \mathbb{D}^{\lambda(p_1)}.$$
In turn applying \cite[Lemma 3.7]{Morse} to $\phi_2$ gives a homotopy equivalence
$$\phi_3: \bar{\mathcal{F}}^{-}\times (0,1)\cup_{\tilde{H}_1^{-}} \mathbb{D}^{\lambda(p_1)}\cup_{\tilde{H}_2^{-}} \mathbb{D}^{\lambda(p_2)}$$
$$\to \bar{\mathcal{F}}^{-}\cup_{\phi_1\circ \tilde{H}_1^{-}} \mathbb{D}^{\lambda(p_1)}\cup_{\phi_2\circ \tilde{H}_2^{-}} \mathbb{D}^{\lambda(p_2)}.$$
Continuing in this manner one eventually obtains a homotopy equivalence between
$$\mathcal{X}_{t,\eta}\sim \bar{\mathcal{F}}^{-}\times (0,1)\cup_{\tilde{H}_1^{-}} \mathbb{D}^{\lambda(p_1)}\cup\dots\cup_{\tilde{H}_m^{-}} \mathbb{D}^{\lambda(p_m)}$$
and 
$$\bar{\mathcal{F}}^{-}\cup_{\phi_1\circ \tilde{H}_1^{-}} \mathbb{D}^{\lambda(p_1)}\cup\dots\cup_{\phi_m\circ \tilde{H}_m^{-}} \mathbb{D}^{\lambda(p_m)}.$$
Then $H_i^{-}:=\phi_i\circ\tilde{H}_i^{-}, i=1,\dots, m$ are the wanted attaching maps. 
\end{enumerate}
\end{proof}

One recovers the formula of Khimshiashvili \cite[Theorem 2.3]{Khimshiashvili} from this Corollary by taking the Euler-Poincaré characteristic. 

\begin{corollary}\label{Euler} With the same assumptions as in Theorem \hyperref[main theorem one]{\ref*{main theorem one} } one has
$$\chi(\bar{\mathcal{F}}_{\eta}^{+})=1-\sum_{i=1}^{m} (-1)^{n+1-\lambda(p_i)},\qquad\chi(\bar{\mathcal{F}}_{\eta}^{-})=1-\sum_{i=1}^{m}(-1)^{\lambda(p_i)}$$
\end{corollary}
\begin{proof} It was established in the Corollary \hyperref[corollary main theorem one]{\ref*{corollary main theorem one} } that
$$\bar{\mathcal{F}}_{t,\eta}^{+}\cup_{H_1^{+}} \mathbb{D}^{n+1-\lambda(p_1)}\cup\dots\cup_{H_m^{+}} \mathbb{D}^{n+1-\lambda(p_m)}$$
is contractible. For any $k=1,\dots, m$ let 
$$\mathcal{M}_k:=\bar{\mathcal{F}}_{t,\eta}^{+}\cup_{H_1^{+}} \mathbb{D}^{n+1-\lambda(p_1)}\cup\dots\cup_{H_k^{+}} \mathbb{D}^{n+1-\lambda(p_k)}.$$
$$\mathcal{M}_0:=\bar{\mathcal{F}}_{t,\eta}^{+}.$$
For any pair $(X, A)$ of topological spaces consider the Euler characteristic 
$$\chi(X,A)=\sum_{n\geq 0} (-1)^n \text{rank } H_n(X, A; \mathbb{Z})$$
in relative singular homology with $\mathbb{Z}$-coefficients. From the long exact sequence of a triple $(X, A, B)$ with  $B\subset A\subset X$ one has that $\chi$ is additive in the sense that 
$$\chi(X, B)=\chi(X, A)+\chi(A, B).$$
Since 
$$\mathcal{M}_0\subset \mathcal{M}_1\subset\dots \subset \mathcal{M}_m$$
is an ascending filtration of topological spaces and since $\mathcal{M}_i\sim \mathcal{M}_{i-1}\cup_{H_i^{+}} \mathbb{D}^{n+1-\lambda(p_i)}$ for $i=1,\dots, m$ one obtains
$$\chi(\mathcal{M}_m, \mathcal{M}_0)=\sum_{i=1}^m \chi(\mathcal{M}_i, \mathcal{M}_{i-1}).$$
Using the Excision Theorem \cite[Theorem 2.20]{Hatcher} applied to the pair $(\mathcal{M}_i, \mathcal{M}_{i-1})$ and the subspace $\mathcal{M}_{i-1}\setminus H_{i-1}^{+}(\partial \mathbb{D}^{n+1-\lambda(p_i)})$ gives 
$$\chi(\mathcal{M}_i, \mathcal{M}_{i-1})=\chi(\mathbb{D}^{n+1-\lambda(p_i)}, \partial\mathbb{D}^{n+1-\lambda(p_i)})=(-1)^{n+1-\lambda(p_i)}$$
(see also \cite[I § 5]{Morse}) hence
$$\chi(\mathcal{M}_m, \mathcal{M}_0)=\sum_{i=1}^{m} (-1)^{n+1-\lambda(p_i)}.$$
One the other hand, applying additivity of $\chi$ to the ascending filtration 
$$\emptyset\subset \mathcal{M}_0=\bar{\mathcal{F}}_{t,\eta}^{+}\subset \mathcal{M}_m$$
of topological spaces one obtains
$$\chi(\bar{\mathcal{F}}_{t,\eta}^{+})=\chi(\mathcal{M}_m)-\chi(\mathcal{M}_m; \mathcal{M}_0).$$
Since $\bar{\mathcal{X}}_{t,\eta}\sim \mathcal{M}_m$ is contractible and since $\bar{\mathcal{F}}_{t,\eta}^{+}\cong \bar{\mathcal{F}}_{\eta}^{+}$ this yields
$$\chi(\bar{\mathcal{F}}_{\eta}^{+})=1-\sum_{i=1}^{m} (-1)^{n+1-\lambda(p_i)}$$
whence the first assertion. The proof of the second assertion is analoguous. 
\end{proof}

\begin{remark} The proof of Khimshiashvili's formula which is known to us\footnote{The original article \cite{khim} is not avaible online} is the one given in \cite[Theorem 2.3]{Khimshiashvili}. It uses Morse theory for manifolds-with-boundaries. Theorem \hyperref[main theorem one]{\ref*{main theorem one} } is not stated in the literature and it does not follow from of the proof \cite[Theorem 2.3]{Khimshiashvili}. 
\end{remark}

Another corollary is the following.

\begin{corollary}\label{contract} In the notation of Theorem \hyperref[main theorem one]{\ref*{main theorem one} } if there exists $t\in V$ with $|t|<t_0$ such that $f_t=F(x,t)$ has no critical points inside the ball of radius $\delta$ centered at the origin then the positive and negative Milnor fibres are contractible. 
\end{corollary}

\begin{remark} If $f: (\mathbb{C}^{n+1},0)\to (\mathbb{C},0)$ is a germ of complex hypersurface singularities then it is never the case that a Milnor fibre (see e.g \cite[Theorem 4.8]{Milnor}) of $f$ is contractible. Indeed the Milnor fibre is contractible if and only the origin is a nonsingular point of $(f^{-1}(0),0)$, according to a theorem \cite[Théorème 3]{Lefschetz} of Norbert A'Campo, generalising a result of Milnor for isolated singularities.
\end{remark}

We give another example of Theorem \hyperref[main theorem one]{\ref*{main theorem one} }.
\begin{example}\label{example stupid} Let $f:\mathbb{R}^2\to\mathbb{R}$ be given by $f(x,y)=y^2-x^3$. Then 
$$F: \mathbb{R}^2\times[0,1]\to\mathbb{R},\qquad F(x,y, t)=y^2-x^3-tx$$
is a $\mathbb{R}$-morsification. Indeed
$$\text{Jac}(f_t)(x,y)=[-(3x^2+t), 2y]$$
has full rank, for any $t\neq 0$. Therefore $f_t: \mathbb{R}^2\to\mathbb{R}$ has no critical points and hence the Milnor fibres $\bar{\mathcal{F}}_{\eta}^{\pm}$ are both contractible by Corollary \hyperref[contract]{\ref*{contract} }. Remark that if we instead consider the $\mathbb{R}$-morsification 
$$G: \mathbb{R}^2\times[0,1]\to\mathbb{R},\qquad G(x,y, t)=y^2-x^3+tx$$ 
then 
$$\text{Jac}(g_t)(x,y)=[t-3x^2, 2y]$$
has for $t\neq 0$ nonmaximal rank in $p_1(t)=(\sqrt{t/3}, 0)$ and in $p_2(t)=(-\sqrt{t/3}, 0)$. Since 
\[
\text{Hess}(g_t)= \begin{pmatrix}
-6x & 0\\
0 & 2 &
\end{pmatrix}
\]
one obtains the Morse indices $\lambda(p_1(t))=1$ and $\lambda(p_2(t))=0$. As a consequence $\bar{\mathcal{F}}_{\eta}^{+}\times(0,1)$ remains contractible after attaching a handle $\mathbb{D}^1\times \mathbb{D}^1$ along $\mathbb{D}^1\times\partial \mathbb{D}^1$ and then attaching a handle $\mathbb{D}^0\times \mathbb{D}^2$  along $\mathbb{D}^0\times \partial \mathbb{D}^2$.
\end{example}

\newpage
\section{On the Homology of The Real Milnor Fibres}\label{cancellation}
According to the previous Theorem \hyperref[main theorem one]{\ref*{main theorem one}} one obtains a contractible space by succesively attaching handles (as in \cite[Notational Definition, p.61]{Goresky}) to (a space having the homotopy type of) the positive (respectively negative) real Milnor fibres. However this gives no information about the manner in which the handles are attached. Moreover it might be the case that there exist $t\in V$ and $1\leq i, j\leq m$ and a pair of handles 
$$(\mathbb{D}^{\lambda(p_i)}\times \mathbb{D}^{n+1-\lambda(p_i)}, \mathbb{D}^{\lambda(p_{j})}\times \mathbb{D}^{n+1-\lambda(p_j)}),$$
such that $\bar{\mathcal{X}}_{t,\eta}$ has the homotopy type of 
$$\bar{\mathcal{F}}_{\eta}^{+}\times (0,1)\cup \left(\bigcup_{1\leq l\neq i, j\leq m}\mathbb{D}^{\lambda(p_l)}\times \mathbb{D}^{n+1-\lambda(p_l)}\right).$$
An example of this phenomenon is provided by Example \hyperref[example stupid]{\ref*{example stupid}} where one knows by using the $\mathbb{R}$-morsification $F$ that the positive Milnor fibre is contractible so that the pair of handles provided by the $\mathbb{R}$-morsification $G$ does not contribute to the homotopy type of $\bar{X}_{t,\eta}$.\\

One therefore concludes that the homology groups of the Milnor fibres cannot be obtained directly from Theorem \hyperref[main theorem one]{\ref*{main theorem one}}. However we shall now state conditions under the validity of which, all of the handles are up to homotopy attached to a real Milnor fibre.\\

Throughout this section we fix Milnor data $(\delta_0,\epsilon_0)$ for $f: (\mathbb{R}^{n+1}, 0)\to (\mathbb{R},0)$.
\subsection{A Lemma}
By restricting the parameter space $U$ one can always find a $\mathbb{R}$-morsification such that the number of nondegenerate critical points lying inside the closed ball of radius $\delta$ is independent of the parameter. 

\begin{lemma}\label{allmen} Let $F: M\times U\to \mathbb{R}$ be a $\mathbb{R}$-morsification of\\
$f:(\mathbb{R}^{n+1}, 0)\to (\mathbb{R}, 0)$. There exists a finite set of nonempty connected semialgebraic subsets 
$$U_1,\dots, U_N\subset U\subset\mathbb{R}^k$$
each containing the origin such that if $t,t'\in U_i$ then $m(t)=m(t')$ (Definition \hyperref[number]{\ref*{number}}) and such that $V\cap U_i=U_i\setminus\{0\}$. That is, the number of critical points of $f_t: M\to\mathbb{R}$ lying inside the ball of radius $\delta$ is independent of $t\in V\cap U_i$ and each critical point is Morse with pairwise distinct critical values.
\end{lemma}
\begin{proof} Since $F: M\times U\to \mathbb{R}$ is a representative of a polynomial map, $M$ and $U$ are semialgebraic. Since
$$\mathcal{P}: t\in U\quad \exists p\in M:\quad \text{rank }\text{Jac}(f_t)(p)=0,\quad \det\text{Hess}(f_t)(p)\neq 0$$
$$\mathcal{Q}: t\in U\quad p,q\in M: \text{rank }\text{Jac}(f_t)(p)=\text{rank }\text{Jac}(f_t)(q)=0,$$
$$\quad p\neq q\implies f_t(p)\neq f_t(q)$$
are semialgebraic conditions one can assume (Definition \hyperref[definition morsification]{\ref*{definition morsification}}) that $V\subset U\setminus\{0\}$ is semialgebraic. Let 
$$\mathscr{C}_V=\{(p,t)\in M\times V\ |\ \text{rank }\text{Jac}(f_t)(p)=0 \}.$$
Then $\mathscr{C}_V$ is a semialgebraic set. Let 
$$\pi: \mathbb{R}^{n+1}\times \mathbb{R}^k\to \mathbb{R}^k,\qquad \pi(x,t)=t$$
denote the standard projection onto the second factor. Then $\pi_{|\mathscr{C}_V}$ is a semialgebraic, continuous map and its fiber dimension is the number of critical points of $f_{t|V}$.  By Hardt's Theorem \cite[Theorem 4.1]{Coste} one can decompose $V$ into a finite union of semialgebraic subsets $V_{1},\dots, V_{L} \subset V$ such that $\pi_{|\mathscr{C}_V}$ is semialgebraically trivial (see e.g \cite[§ 4.1.1 ]{Coste}) over each $V_{i}$. The fiber dimension $\dim \pi_{|\mathscr{C}}^{-1}(t)$ is therefore independent of $t\in V_i$. In other words the number of critical points of $f_t$ is independent of $t\in V_i$, each critical point is Morse and the critical values are pairwise distinct. Since $(0,\dots, 0)\in\bar{V}$ by the Morse Lemma \cite[Lemma 2.2]{Morse} there exists a subcollection $\{V_1,\dots, V_N\}\subset \{V_{1},\dots, V_L\}$ such that $(0,\dots, 0)\in\bar{V}_i, i=1,\dots, N$.  Every semialgebraic set has a decomposition into a finite number of connected semialgebraic components by \cite[Theorem 2.4.5]{roy} so we can assume that $V_1,\dots, V_N$ are connected. Then $U_i=V_i\cup\{0\}, i=1,\dots,N$ satisfy the conditions of the Lemma.
\end{proof}

\begin{remark}
One can thus always obtain a $\mathbb{R}$-morsification $F$ such that that $f_t: M \to \mathbb{R}$ is a Morse function with exactly $m=m(t)$ critical points with distinct critical values for \emph{all} $t\in U\setminus\{0\}$. Indeed this follows from the previous Lemma \hyperref[allmen]{\ref*{allmen}} since one can replace a given $\mathbb{R}$-morsification $F: M\times U\to\mathbb{R}$ by a new $\mathbb{R}$-morsification $F': M\times U_i\to \mathbb{R}$ where $U_i$ is one of the connected semialgebraic sets given by Lemma \hyperref[allmen]{\ref*{allmen}} and where the dense set $V\subset U$ is replaced by by $U_i\setminus \{0\}$.
\end{remark}

\subsection{The Local Milnor Fibres}\label{subsection local} We shall now introduce the local Milnor fibres. In this subsection $F: M\times U\to \mathbb{R}$ is a $\mathbb{R}$-morsification of $f$ such that $m(t)$ is independent of $t$, denoted $m=m(t)$, and such that $V=U\setminus\{0\}$.\\

Let $t\in V$. By the Morse Lemma \cite[Lemma 2.2]{Morse} each of the critical points $p_i=p_i(t)\in\bar{\mathbb{B}}_{\delta}$ of $f_t: M\to \mathbb{R}$ defines a quadratic singularity in the fibre $f_t^{-1}(s_i)$, where $s_i=f_t(p_i)$. One can therefore apply the Fibration Theorem of Milnor \hyperref[Milnor fibrations]{\ref*{Milnor fibrations}} \cite[Theorem 4.2]{Seade} to deduce the following.\\

For each $i=1,\dots, m$ there exists a $\delta_{i,0}=\delta_{i,0}(t)>0$ such that for any $\delta_i\in (0,\delta_{i,0}]$ there exist an $\epsilon_{i,0}>0$ such that for any $\epsilon_i\in (0, \epsilon_{i,0})$ the restrictions 
\begin{equation}
\label{one}
f_t: f_t^{-1}((s_i, s_i+ \epsilon_i])\cap \mathbb{B}_{\delta_i}(p_i)\to (s_i,s_i+\epsilon_i]
\end{equation}
and 
\begin{equation}
\label{two}
f_t: f_t^{-1}([s_i-\epsilon_i, s_i))\cap \mathbb{B}_{\delta_i}(p_i)\to [s_i-\epsilon_i,s_i)
\end{equation}
are projections of trivial (topological) fibrations where $\mathbb{B}_{\delta_i}(p_i)$ is the open ball centered at $p_i$ and of radius $\delta_i$. 
The pair $(\delta_{i,0}, \epsilon_{i,0})=(\delta_{i,0}(t), \epsilon_{i,0}(t))$ will be called \emph{local Milnor data} at the point $p_i$.\\

\begin{remark} We shall write $\mathbb{B}_{\delta_i}$ instead of $\mathbb{B}_{\delta_i}(p_i)$ when no confusion is possible.
\end{remark}

\begin{definition}\label{local} Let $t\in V$ be fixed and let $\delta_i\in (0,\delta_{i,0}]$ and $\eta_i\in (0,\epsilon_{i,0}]$. The fibers
$$\mathcal{F}_{i,loc}^{+}:=f_t^{-1}(s_i+\eta_i)\cap \mathbb{B}_{\delta_i},\qquad i=1,\dots,m$$
$$\mathcal{F}_{i,loc}^{-}:=f_t^{-1}(s_i-\eta_i)\cap \mathbb{B}_{\delta_i},\qquad i=1,\dots,m$$
of the fibrations \hyperref[one]{\ref*{one}} and \hyperref[two]{\ref*{two}} are called the local positive (respectively negative) open Milnor fibres at the point $p_i=p_i(t)$ with respect to the $\mathbb{R}$-morsification $F$. Their closures
$$\bar{\mathcal{F}}_{i,loc}^{\pm}:=f_t^{-1}(s_i+\eta_i)\cap \bar{\mathbb{B}}_{\delta_i},\qquad i=1,\dots,m$$
$$\bar{\mathcal{F}}_{i,loc}^{\pm}:=f_t^{-1}(s_i-\eta_i)\cap \bar{\mathbb{B}}_{\delta_i},\qquad i=1,\dots,m$$
are called the local positive (respectively negative) closed Milnor fibres at the point $p_i$ with respect to $F$.
\end{definition}

\begin{figure}[h]
\begin{tikzpicture}[scale=0.5]

\filldraw[gray, opacity=0.3] (0,0) circle [radius=4.5 cm];
\draw (0,0) circle [radius=4.5 cm];
\draw (0,-4.5) .. controls (-0.5, 0) .. (0, 4.5);

\draw (2,0) .. controls (2,0.5) and (1,2) .. (1.4,4.3);
\draw (2,0) .. controls (2,-0.5) and (1,-2) .. (1.4,-4.3);

\draw (-2,0) .. controls (-1,1) and (-1,2) .. (-1.4,4.3);
\draw (-2,0) .. controls (-1,-1) and (-1,-2) .. (-1.4,-4.3);

\draw (2.5,-3.75) ..controls (3.5, 0) .. (2.5, 3.75);
\draw (-2.5,-3.75) .. controls (-3.5,0) .. (-2.5,3.75);

\draw (-3,-4.5) node [below right]{$\bar{\mathbb{B}}_{\delta}$}; 

\filldraw[gray, opacity=0.6] (2,0) circle [radius=0.8 cm];
\filldraw[gray, opacity=0.6] (-2,0) circle [radius=0.8];
\draw (2,0) circle [radius=0.8 cm];
\draw (-2,0) circle [radius=0.8 cm];

\draw (2.2, -1) node[below right] {$\bar{\mathbb{B}}_{\delta_j}$};
\draw (-2.2, -1) node[below right] {$\bar{\mathbb{B}}_{\delta_i}$};

\filldraw [black] (-2,0) circle (2pt) node[right] {$+$};
\filldraw [black] (-2,0) circle (2pt) node[left] {$-$};
\filldraw [black] (2,0) circle (2pt) node[right] {$+$};
\filldraw [black] (2,0) circle (2pt) node[left] {$-$};

\begin{flushbottom}

\draw[->] (0, -6) --  (0,-8) node[above right] {$f_{t}$};



\draw (-3,-10) -- (3,-10);

\filldraw [black] (0,-10) circle (2pt) node[below] {$0$};
\filldraw [black] (-3,-10) circle (2pt) node[below] {$-\eta$};
\filldraw [black] (3,-10) circle (2pt) node[below] {$\eta$};
\filldraw [black] (-1.5,-10) circle (2pt) node[below] {$s_i$};
\filldraw [black] (1.5,-10) circle (2pt) node[below] {$s_j$};

\end{flushbottom}
\end{tikzpicture}
\caption{A set of two local Milnor fibres}\label{fig}
\end{figure}
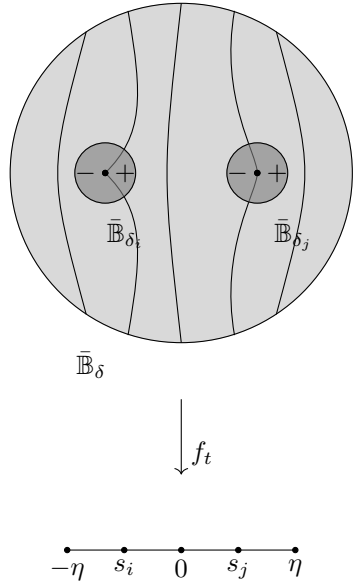

The following schematic Figure \hyperref[fig]{\ref*{fig}} serves to illustrate the local Milnor fibres and also the proof of the next Theorem \hyperref[main theorem two]{\ref*{main theorem two}}.

\subsection{The Main Theorem}
The main result of this section is the following theorem. We assume throughout that $\bar{\mathcal{F}}_{\eta}^{+}$ and $\bar{\mathcal{F}}_{\eta}^{-}$ are nonempty

\begin{theorem}\label{main theorem two} Let $F: M\times U\to \mathbb{R}$ be a $\mathbb{R}$-morsification of a germ of isolated singularity $f: (\mathbb{R}^{n+1}, 0)\to (\mathbb{R},0)$ with $n\equiv 1\quad(\text{mod }2)$ and $n>1$. Let $\delta\in (0,\delta_0]$ and $\eta\in (0,\epsilon_0]$ and recall that $m=m(t)$ for all $t\in V$ (where $m(t)$ is as in Definition \hyperref[number]{\ref*{number}} and $V$ is as in Definition \hyperref[definition morsification]{\ref*{definition morsification}}). Let $t_0$ be as in Lemma \hyperref[main lemma one]{\ref*{main lemma one}} and let $t\in V\cap\{|t|\leq t_0\}$ be fixed. Denote by $p_1,\dots, p_m\in\bar{\mathbb{B}}_{\delta}$ the critical points of $f_t: M\to \mathbb{R}$ lying inside the ball of radius $\delta$ and let $\lambda(p_i)$ denote their Morse indices. If $\lambda(p_i)=(n+1)/2,i=1,\dots, m$ then 
$$H_k(\bar{\mathcal{F}}_{\eta}^{+})=H_k(\bar{\mathcal{F}}_{\eta}^{-})\cong H_k(\bigvee_{i=1}^m \mathbb{S}^{(n-1)/2}),\qquad k\geq 0.$$
\end{theorem} 
\begin{proof}
Using the Isotopy Lemma \hyperref[Thom]{\ref*{Thom}} (\cite[Proposition 11.1]{Mather}) together with Corollary \hyperref[corollary main theorem one]{\ref*{corollary main theorem one}} we first show that the contractible space $\bar{\mathcal{X}}_{t,\eta}$ is obtained by attaching discs $\mathbb{D}^{\lambda(p_i)}$ along their boundaries $\partial \mathbb{D}^{\lambda(p_i)}$ to certain spaces which by using \cite[Proposition 3.2]{Goresky} are shown to be homotopy equivalent to different copies of the negative Milnor fibre of $f$. Using the Excision Theorem \cite[Theorem 2.20]{Hatcher} one is the able to obtain the integral homology groups of the negative Milnor fibre. The case with the positive Milnor fibres is then obtained by considering $-f_t$ instead of $f_t$. 
\begin{enumerate}
\item\label{step one thing} Fix $t\in V$. Let $\delta_i=\delta_i(t)\in (0,\delta_{i,0}(t))$ and $\epsilon_i=\epsilon_i(t)\in (0,\epsilon_{i,0}(\delta_i))$ be local Milnor data as in Definition \hyperref[local]{\ref*{local}}. For each $i=1,\dots, m$ one can assume by the Morse Lemma \cite[Lemma 2.2]{Morse} that $\delta_i$ is chosen such that there exists local coordinates $x_1,\dots, x_{n+1}: \mathbb{B}_{\delta_i}\to \mathbb{R}^{n+1}$ with $x_j(p_i)=(0,\dots, 0), j=1,\dots, n+1$ such that
$$f_{|\mathbb{B}_{\delta_i}}=s_i-\sum_{j=1}^{\lambda} x_j^{2}+\sum_{l=\lambda+1}^{n+1} x_i^2,$$
up to precomposing $f$ with a diffeomorphism, where $\lambda:=\lambda(p_i), i=1,\dots, m$. Put 
$$\mathcal{Z}_{i}=f_t^{-1}(s_i-\epsilon_i, s_i+\epsilon_i)\cap \bar{\mathbb{B}}_{\delta},$$
$$\bar{\mathcal{Z}}_{i}=f_t^{-1}([s_i-\epsilon_i, s_i+\epsilon_i])\cap \bar{\mathbb{B}}_{\delta},$$
$$\mathcal{Y}_{i}=\mathcal{Z}_i\cap\mathbb{B}_{\delta_i}=f_t^{-1}(s_i-\epsilon_i, s_i+\epsilon_i)\cap \mathbb{B}_{\delta_i}.$$
For each $i=1,\dots, m$ choose real positive numbers $\gamma_i>0$ such that the interval $[s_i-\gamma_i, s_i+\gamma_i]$ contains no critical value of $f_{t|\bar{\mathbb{B}}_{\delta}}$ except $s_i$ and such that 
$$s_i+\gamma_i=s_{i+1}-\gamma_{i+1},\qquad i=1,\dots, m-1,$$
$$-\eta=s_1-\gamma_1,\quad \eta=s_m+\gamma_m$$
Put 
$$\mathcal{X}_{i}=f_t^{-1}(s_i-\gamma_i, s_i+\gamma_i)\cap \bar{\mathbb{B}}_{\delta},$$
$$\bar{\mathcal{X}}_{i}=f_t^{-1}([s_i-\gamma_i, s_i+\gamma_i])\cap \bar{\mathbb{B}}_{\delta}$$
Then the union of the sets $\bar{\mathcal{X}}_i, i=1,\dots, m$ is the manifold with corners $\bar{\mathcal{X}}_{\eta}:=\bar{\mathcal{X}}_{t,\eta}$ of Lemma \hyperref[main lemma one]{\ref*{main lemma one}}. In particular it is contractible. By Lemma \hyperref[manifold lemma]{\ref*{manifold lemma}} one has that each $\bar{\mathcal{X}}_{i}$ and $\bar{\mathcal{Z}}_i$ are smooth manifolds with corners. By smoothening their corners one deduces that $\bar{\mathcal{X}}_i\sim \mathcal{X}_i$ and $\bar{\mathcal{Z}}_i\sim \mathcal{Z}_i$. Let $i\in\{1,\dots, m\}$. We shall now use Goresky and MacPherson's method of ``moving the wall'' to deduce that $\bar{\mathcal{X}}_i$ is homeomorphic to $\bar{\mathcal{Z}}_i$. For this consider 
$$\mathcal{W}=\{(x,r)\in \mathbb{R}\times [0, 1]\ |\ x\in [s_i-r\gamma_i-(1-r)\epsilon_i, s_i+r\gamma_i+(1-r)\epsilon_i]\}$$
and endow $\mathcal{W}$ with its natural Whitney stratification. If $\pi: \mathbb{R}\times\mathbb{R}\to \mathbb{R}$ is the standard projection $\pi(x,r)=r$ then $\pi_{|\mathcal{W}}$ is a proper submersion by construction. Now, by the choice of $t_0\in\mathbb{R}$ one has that $f_t: \bar{\mathbb{B}}_{\delta}\to \mathbb{R}$ is transverse to $[-\epsilon, \epsilon]$ hence $f_{t|\bar{\mathbb{B}}_{\delta}}$ is transverse to  
$$\mathcal{W}_r=[s_i-r\gamma_i-(1-r)\epsilon_i, s_i+r\gamma_i+(1-r)\epsilon_i]\subset [-\epsilon, \epsilon]$$
for all $r\in [0,1]$. One can therefore apply \cite[Theorem 4.8]{Goresky}. This gives a stratum-preserving homeomorphism
$$f_t^{-1}([s_i-\epsilon_i, s_i+\epsilon_i])\cap \bar{\mathbb{B}}_{\delta}\cong f_t^{-1}([s_i-\gamma_i, s_i+\gamma_i])\cap \bar{\mathbb{B}}_{\delta}$$
As a consequence  $\bar{\mathcal{X}}_i\cong \bar{\mathcal{Z}}_i$ and therefore also $\mathcal{X}_i\sim \mathcal{Z}_i$.

\item\label{step two thing} For any $i=1,\dots, m$ consider the map
$$f_t: \mathcal{Z}_{i}\setminus \mathcal{Y}_{i}\to (s_i-\epsilon_i,s_i+\epsilon_i).$$
It is proper since its fibers are of the form
$$f_{t| \mathcal{Z}_{i}\setminus\mathcal{Y}_{i}}^{-1}(s)=f_t^{-1}(s)\cap (\bar{\mathbb{B}}_{\delta}\setminus \mathbb{B}_{\delta_i})$$
and hence are compact. From the assumption 
$$\lambda(p_1)=\dots=\lambda(p_m)=\frac{n+1}{2}$$
follows that if we choose $\eta_i\in (0,\epsilon_i], i=1,\dots, m$ then local Milnor fibres (Definition \hyperref[local]{\ref*{local}}) are of the form
$$\bar{\mathcal{F}}_{i,loc}^{+}\cong \{x\in\bar{\mathbb{B}}_{\delta_i}\ |\ -\sum_{j=1}^{(n+1)/2} x_j^2+\sum_{j=(n+3)/2}^{n+1} x_j^2=\eta_i\}\cong \bar{\mathcal{F}}_{i,loc}^{-}$$
In particular the boundaries of the local Milnor fibres $\bar{\mathcal{F}}_{t, i,loc}^{\pm}$ are nonempty so the fibers of $f_{t| \mathcal{Z}_{i}\setminus\mathcal{Y}_{i}}$ are nonempty. Note that the only stratified critical point of $f_{t| \mathcal{Z}_{i}}$ is by construction $p_i\in\mathcal{Y}_{i}$. Therefore, in order to show that $f_{t| \mathcal{Z}_{i}\setminus\mathcal{Y}_{i}}$ is a proper stratified submersion it suffices to show that 
$$f_t: f_t^{-1}(s_i-\epsilon_i, s_i+\epsilon_i)\cap\mathbb{S}_{\delta_i}\to \mathbb{R}$$
is a submersion, that is, that the map
$$f_t: \mathbb{S}_{\delta_i}\to\mathbb{R}$$ 
has no critical values in $(s_i-\epsilon_i, s_i+\epsilon_i)$. It has no critical values $s\neq s_i$ because 
$$f_t: f_t^{-1}((s_i, s_i+\epsilon_i])\cap\bar{\mathbb{B}}_{\delta_i}\to (s_i, s_i+\epsilon_i],$$
$$f_t: f_t^{-1}((s_i-\epsilon_i, s_i])\cap\bar{\mathbb{B}}_{\delta_i}\to [s_i-\epsilon_i, s_i)$$
are stratified submersions and because $f_t^{-1}(s_i)\cap\mathbb{S}_{\delta_i}$ is transverse by construction (Definition \hyperref[local]{\ref*{local}}) and so $f_t: \mathbb{S}_{\delta_i(t)}\to \mathbb{R}$ has no critical point above $s_i$. All said $f_{t| \mathcal{Z}_{i}\setminus\mathcal{Y}_{i}}$ is a proper stratified submersion so by Thom's First Isotopy Lemma \hyperref[Thom]{\ref*{Thom}} (\cite[Proposition 11.1]{Mather}),
$$\mathcal{Z}_{i}\setminus \mathcal{Y}_{i}\cong f_t^{-1}(s)\cap (\bar{\mathbb{B}}_{\delta}\setminus \mathbb{B}_{\delta_i})\times (s_i-\epsilon_i, s_i+\epsilon_i),$$
for all $s\in (s_i-\epsilon_i, s_i+\epsilon_i)$. This yields in particular a deformation retraction
$$\mathcal{Z}_{i}\sim f_t^{-1}(s)\cap (\bar{\mathbb{B}}_{\delta}\setminus \mathbb{B}_{\delta_i})\cup \mathcal{Y}_{i}=$$
$$=f_t^{-1}(s)\cap \bar{\mathbb{B}}_{\delta} \cup \mathcal{Y}_{i},\qquad s\in (s_i-\epsilon_i, s_i+\epsilon_i).$$
\item\label{step tve thing}  We claim that $\mathcal{Z}_i\sim \mathcal{Z}_{i-1}$ for all $i=1,\dots, m$. To prove this we shall show that both $\mathcal{Z}_i$ and $\mathcal{Z}_{i-1}$ deformation retracts to homeomorphic spaces. From the Morse Lemma \cite[Lemma 2.2]{Morse} follows that there exists a diffeomorphism $\mathcal{Y}_{i-1}\cong \mathcal{Y}_i$. By the previous step \hyperref[step two thing]{\ref*{step two thing}} it suffices therefore to show that
$$f_t^{-1}(s_i-\eta_i)\cap (\bar{\mathbb{B}}_{\delta}\setminus \mathbb{B}_{\delta_i})\cong f_t^{-1}(s_{i-1}+\eta_{i-1})\cap (\bar{\mathbb{B}}_{\delta}\setminus \mathbb{B}_{\delta_{i-1}}).$$
The interval $[s_i-\gamma_i, s_i)$ contains no critical values of the proper map $f_{t|\bar{\mathbb{B}}_{\delta}}$. There exists therefore by \cite[Proposition 3.2]{Goresky} a stratum-preserving\footnote{By the proof of \cite[Proposition 3.2]{Goresky}, which is a simple consequence of \cite[Theorem 4.4]{Goresky}, it follows that $\phi_i$ is in fact smooth on strata. This follows also from \cite[Remark 7.2]{Goresky} where one uses the First Isotopy Lemma to lift a certain smooth vector field to construct the homeomorphism $\phi_i$.}  homeomorphism
$$\phi_i: f_t^{-1}((-\infty, s_{i-1}+\eta_{i-1}])\cap\bar{\mathbb{B}}_{\delta}\to f_t^{-1}((-\infty, s_i-\eta_i])\cap\bar{\mathbb{B}}_{\delta}.$$
which hence restricts to a homeomorphism
$$\phi_{i|A_i}: A_i:=f_t^{-1}(s_{i-1}+\eta_{i-1})\cap\bar{\mathbb{B}}_{\delta}\to B_i=f_t^{-1}(s_i-\eta_i)\cap\bar{\mathbb{B}}_{\delta}$$
which preserves strata. Since $\mathcal{F}_{i-1, loc}^{+}\subset\text{int}(A_i)$ it follows that $\phi_{i|\mathcal{F}_{i-1, loc}^{+}}$ is a homeomorphism onto its image. Since $\mathcal{F}_{i-1, loc}^{+}$ and $\mathcal{F}_{i,loc}^{-}$ are diffeomorphic by \cite[Lemma 2.2]{Morse} it follows that the image of $\phi_{i|\mathcal{F}_{i-1, loc}^{+}}$ is homeomorphic to $\mathcal{F}_{i, loc}^{-}$. Then
$$\tilde{\phi}_{i |A_i\setminus \mathcal{F}_{i-1,loc}^+}: A_i\setminus \mathcal{F}_{i-1,loc}^+\to B_i\setminus \phi_i (\mathcal{F}_{i-1, loc}^+)\cong B_i\setminus \mathcal{F}_{i,loc}^-$$
gives the wanted homeomorphism. This yields
$$\mathcal{Z}_{i}\sim f_t^{-1}(s_i-\eta_i)\cap (\bar{\mathbb{B}}_{\delta}\setminus \mathbb{B}_{\delta_i})\cup \mathcal{Y}_{i}$$
$$\cong f_t^{-1}(s_{i-1}+\eta_{i-1})\cap (\bar{\mathbb{B}}_{\delta}\setminus \mathbb{B}_{\delta_{i-1}})\cup \mathcal{Y}_{i-1}\sim \mathcal{Z}_{i-1}.$$
\item\label{step three thing} By the proof of Theorem \hyperref[main theorem one]{\ref*{main theorem one}} and Corollary \hyperref[corollary main theorem one]{\ref*{corollary main theorem one}} applied to 
$$f_t: f_t^{-1}([s_i-\epsilon_i, s_i+\epsilon_i])\cap\bar{\mathbb{B}}_{\delta}\to \mathbb{R}$$
there exists an attaching map 
$$h_i^{-}: \partial\mathbb{D}^{\lambda}\to f_t^{-1}(s_i-\eta_i)\cap\bar{\mathbb{B}}_{\delta_i}$$
and a homotopy equivalence 
$$\mathcal{Y}_{i}\sim f_t^{-1}(s_i-\eta_i)\cap\bar{\mathbb{B}}_{\delta_i}\cup_{h_i^{-}} \mathbb{D}^{\lambda}$$
where one uses the fact that $\mathcal{Y}_i\sim f_t^{-1}(s_i-\epsilon_i, s_i+\epsilon_i)\cap \bar{\mathbb{B}}_{\delta}$ which follows from the Tubular Neighborhood Theorem \cite[Proposition 3.42]{Hatcher} for topological manifolds with boundaries. Applying instead Theorem \hyperref[main theorem one]{\ref*{main theorem one}} and Corollary \hyperref[corollary main theorem one]{\ref*{corollary main theorem one}} to 
$$-f_t: f_t^{-1}([s_i-\epsilon_i, s_i+\epsilon_i])\cap\bar{\mathbb{B}}_{\delta}\to \mathbb{R}$$
one obtains an attaching map
$$h_i^{+}: \partial\mathbb{D}^{\lambda}\to f_t^{-1}(s_i+\eta_i)\cap\bar{\mathbb{B}}_{\delta_i}$$
and a homotopy equivalence
$$\mathcal{Y}_{i}\sim f_t^{-1}(s_i+\eta_i)\cap\bar{\mathbb{B}}_{\delta_i}\cup_{h_i^{+}} \mathbb{D}^{\lambda}$$ 
Applying the previous step \hyperref[step two thing]{\ref*{step two thing}} yields deformation retractions
\[\mathcal{Z}_i\sim  \left\{
\begin{array}{ll}
      f_t^{-1}(s_i+\eta_i)\cap\bar{\mathbb{B}}_{\delta}\cup_{h_i^{+}}  \mathbb{D}^{\lambda},\\
      f_t^{-1}(s_i-\eta_i)\cap\bar{\mathbb{B}}_{\delta}\cup_{h_i^{-}} \mathbb{D}^{\lambda}.
\end{array} \right. \]

\item\label{step four thing} Using \cite[Proposition 3.2]{Goresky} one obtains homeomorphisms
$$\phi_i^{+}: f_t^{-1}(s_i+\eta_i)\cap \bar{\mathbb{B}}_{\delta}\to f_t^{-1}(s_i+\gamma_i)\cap \bar{\mathbb{B}}_{\delta}$$
$$\phi_i^{-}: f_t^{-1}(s_i-\eta_i)\cap \bar{\mathbb{B}}_{\delta}\to f_t^{-1}(s_i-\gamma_i)\cap \bar{\mathbb{B}}_{\delta}$$
so that applying \cite[Lemma 3.7]{Morse} and step \hyperref[step three thing]{\ref*{step three thing}} one obtains homotopy equivalences 
\[\mathcal{Z}_i\sim  \left\{
\begin{array}{ll}
      f_t^{-1}(s_i+\gamma_i)\cap\bar{\mathbb{B}}_{\delta}\cup_{\phi_i^{+}\circ h_i^{+}} \mathbb{D}^{\lambda},\\
      f_t^{-1}(s_i-\gamma_i)\cap\bar{\mathbb{B}}_{\delta}\cup_{\phi_i^{-}\circ h_i^{-}} \mathbb{D}^{\lambda}.
\end{array} \right. \]

\item\label{fick} By the first step \hyperref[step one thing]{\ref*{step one thing}} one has $\bar{\mathcal{X}}_i\sim \mathcal{X}_i$ and $\mathcal{X}_i\sim \mathcal{Z}_i$ for all $i=1,\dots, m$. By step \hyperref[step tve thing]{\ref*{step tve thing}} one has $\mathcal{Z}_i\sim\mathcal{Z}_{i-1}$ for all $i=2,\dots, m$. Hence 
$$\bar{\mathcal{X}}_i\sim f_t^{-1}(-\eta)\cap\bar{\mathbb{B}}_{\delta}\cup_{\tilde{H}_i} \mathbb{D}^{\lambda}$$
for certain attaching maps $\tilde{H}_i$, by step \hyperref[step four thing]{\ref*{step four thing}}.
\item\label{last} Consider the long exact sequence of the pair $(\bar{\mathcal{X}}_{t,\eta},\bar{\mathcal{F}}_{t,\eta}^{-})$ in integral homology:
$$\to\tilde{H}_n(\bar{\mathcal{F}}_{t,\eta}^{-})\to \tilde{H}_n(\bar{\mathcal{X}}_{t,\eta})\to H_{n}(\bar{\mathcal{X}}_{t,\eta},\bar{\mathcal{F}}_{t,\eta}^{-})$$
Since by Lemma \hyperref[main lemma one]{\ref*{main lemma one}} $\bar{\mathcal{X}}_{t,\eta}$ is contractible this gives
$$\tilde{H}_k(\bar{\mathcal{F}}_{t,\eta}^{-})=H_{k+1}(\bar{\mathcal{X}}_{t,\eta}, \bar{\mathcal{F}}_{t,\eta}^{-}),\qquad\forall k\geq 0$$
By the previous step \hyperref[fick]{\ref*{fick}} and the fact that $\bar{\mathcal{X}}_i\sim \mathcal{X}_i$ gives 
$$\bar{\mathcal{X}}_{\eta}\sim \bar{\mathcal{F}}_{t, \eta}^{-}\cup_{\tilde{H}_1} \mathbb{D}^{\lambda}\cup\dots\cup_{\tilde{H}_m} \mathbb{D}^{\lambda}$$
where
$$\tilde{H}_i: \partial\mathbb{D}^{\lambda}\to \bar{\mathcal{F}}_{t, \eta}^{-},\qquad i=1,\dots,m$$
are attaching maps. Therefore
$$\tilde{H}_k(\bar{\mathcal{F}}_{t,\eta}^{-})\cong H_{k+1}(\bar{\mathcal{F}}_{t, \eta}^{-}\cup_{\tilde{H}_1} \mathbb{D}^{\lambda}\cup\dots\cup_{\tilde{H}_m} \mathbb{D}^{\lambda}, \bar{\mathcal{F}}_{t,\eta}^{-}).$$
Since $\lambda>0$ the closure of 
$$\bar{\mathcal{F}}_{\eta}^{-}\setminus \bigcup_{i=1}^m \tilde{H}_i(\partial\mathbb{D}^{\lambda})$$
is contained in the interior of $\bar{\mathcal{F}}_{\eta}^{-}$ so using the Excision Theorem \cite[Theorem 2.20]{Hatcher} yields
$$\tilde{H}_k(\bar{\mathcal{F}}_{t,\eta}^{-})=\bigcup_{i=1}^m H_{k+1}(\mathbb{D}^{\lambda}, \partial \mathbb{D}^{\lambda})$$
$$=\bigoplus_{i=1}^m \mathbb{Z},\qquad k=\lambda-1.$$
Therefore 
$$H_k(\bar{\mathcal{F}}_{\eta}^{-})=H_k(\bigvee_{i=1}^m \mathbb{S}^{\lambda-1})$$
where we have used the fact that $\bar{\mathcal{F}}_{t,\eta}^{-}\cong \bar{\mathcal{F}}_{\eta}^{-}$ by Lemma \hyperref[main lemma one]{\ref*{main lemma one}}.
\end{enumerate}
\end{proof}

\subsection{Consequences}
We now discuss some consequences.

\begin{corollary} Keep the assumptions of the Theorem \hyperref[main theorem two]{\ref*{main theorem two}}. The Poincaré polynomial of the real Milnor fibres are
$$\beta(\bar{\mathcal{F}}_{\eta}^{+})=1+m u^{n-\lambda}$$
and 
$$\beta(\bar{\mathcal{F}}_{\eta}^{-})=1+m u^{\lambda-1}$$
\end{corollary}

The following result, which follows from \hyperref[main theorem one]{\ref*{main theorem one}}, is of considerable use when one considers e.g. $ADE$-singularities.

\begin{corollary}\label{quadratic case} Keep the notations of Theorem \hyperref[main theorem one]{\ref*{main theorem one}} and suppose $f_t: M\to\mathbb{R}$ has a unique critical point $p$ of index $\lambda(p)$. If $\bar{\mathcal{F}}_{\eta}^{+}$ is nonempty then 
\[\beta(\bar{\mathcal{F}}_{\eta}^{+})=\left\{
\begin{array}{ll}
      1+u^{n-\lambda(p)},\quad \text{if }\lambda(p)\leq n,\\
      1,\quad\text{if } \lambda(p)=n+1.
\end{array} \right. \]
If $\bar{\mathcal{F}}_{\eta}^{-}$ is nonempty then
\[\beta(\bar{\mathcal{F}}_{\eta}^{-})=\left\{
\begin{array}{ll}
      1+u^{\lambda(p)-1},\quad\text{if } \lambda(p)\geq 1,\\
      1,\quad\text{if } \lambda(p)=0.
\end{array} \right. \]
\end{corollary}
\begin{proof} We shall treat the second assertion, the first assertion follows by considering $-f_t$ instead of $f_t$ and repeating the arguments below. By the Corollary \hyperref[corollary main theorem one]{\ref*{corollary main theorem one}} there is an attaching map 
$$h: \partial\mathbb{D}^{\lambda(p)}\to \bar{\mathcal{F}}_{t,\eta}^{-}$$
such that 
$$\bar{\mathcal{X}}_{t,\eta}\sim \bar{\mathcal{F}}_{t, \eta}^{-}\cup_{\partial\mathbb{D}^{\lambda(p)}} \mathbb{D}^{\lambda(p)}$$
is contractible. If $\lambda(p)=0$ then $\mathbb{D}^{\lambda}$ is a point attached to $\bar{\mathcal{F}}_{t,\eta}^{-}$ along the empty set. Then either $\bar{\mathcal{F}}_{t,\eta}^{-}=\emptyset$ in which case $\bar{\mathcal{F}}_{\eta}^{-}=\emptyset$, or $\bar{\mathcal{F}}_{t,\eta}^{-}$ is contractible in which case $\beta(\bar{\mathcal{F}}_{\eta}^{-})=1$. So suppose $\bar{\mathcal{F}}_{\eta}^{-}\neq\emptyset$ and that $\lambda(p)>0$. Then $\bar{\mathcal{F}}_{t, \eta}^{-}$ is nonempty as well so one can apply the long exact sequence in homology \cite[Proposition 4.12]{Ebeling}, yielding 
$$\tilde{H}_k(\bar{\mathcal{F}}_{\eta}^{-})\cong H_{k+1}(\bar{\mathcal{X}}_{t,\eta}, \bar{\mathcal{F}}_{t,\eta}^{-}).$$
Since $\lambda(p)>0$ the closure of $\bar{\mathcal{F}}_{t,\eta}^{-}\setminus h(\partial \mathbb{D}^{\lambda(p)})$ is contained in the interior of 
$\bar{\mathcal{F}}_{t,\eta}^{-}$ so one can use the Excision Theorem \cite[Theorem 2.20]{Hatcher}\footnote{One can also take a neighbourhood $A\supset \bar{\mathcal{F}}_{t,\eta}^{-}$ such that $\text{int}(A)\cup \mathbb{D}^{\lambda(p)}$ deformation retracts to $\bar{\mathcal{F}}_{t, \eta}^{-}\cup_{\partial\mathbb{D}^{\lambda(p)}} \mathbb{D}^{\lambda(p)}$
and such that $A\cap \mathbb{D}^{\lambda(p)}=\partial\mathbb{D}^{\lambda(p)}$ and apply the Excision Theorem to $A$ and $B=\mathbb{D}^{\lambda(p)}$.}. This gives
$$\tilde{H}_k(\bar{\mathcal{F}}_{\eta}^{-})=H_{k+1}(\bar{\mathcal{F}}_{t, \eta}^{-}\cup_{\partial\mathbb{D}^{\lambda(p)}} \mathbb{D}^{\lambda(p)}, \bar{\mathcal{F}}_{t,\eta}^{-})$$
$$\cong H_{k+1}(\mathbb{D}^{\lambda(p)}, \partial\mathbb{D}^{\lambda(p)}).$$
The claim follows.
\end{proof}

\subsection{Real Vanishing Cycles}
In this subsection $F: M\times U\to \mathbb{R}$ is a $\mathbb{R}$-morsification such that $m=m(t)$ and such that $V=U\setminus\{0\}$, as in subsection \hyperref[subsection local]{\ref*{subsection local}}. We assume furthermore that $\bar{\mathcal{F}}_{\eta}^{+}$ and $\bar{\mathcal{F}}^{-}$ are nonempty. Fix $t\in V$ and consider for each $i=1,\dots, m$ the local Milnor fibres (Definition \hyperref[local]{\ref*{local}})
$$\mathcal{F}_{i,loc}^{+}=f_t^{-1}(s_i+ \eta_i)\cap \mathbb{B}_{\delta_i}$$
$$\mathcal{F}_{i,loc}^{-}=f_t^{-1}(s_i-\eta_i)\cap \mathbb{B}_{\delta_i}$$
at the critical point $p_i=p_i(t)\in \bar{\mathbb{B}}_{\delta}$ of $f_t=F(\cdot, t)$. Let $\lambda(p_i)$ be the Morse index of $f_t$ at $p_i$. Then Corollary \hyperref[quadratic case]{\ref*{quadratic case}} gives 
\[H_k(\mathcal{F}_{i,loc}^{+})=\left\{
\begin{array}{ll}
      \mathbb{Z},\quad k=0, n-\lambda(p_i),\\
      \{0\},\quad \text{otherwise}
\end{array} \right. \]

\[H_k(\mathcal{F}_{i,loc}^{-})=\left\{
\begin{array}{ll}
      \mathbb{Z},\quad k=0, \lambda(p_i)-1,\\
      \{0\},\quad \text{otherwise}
\end{array} \right. \]
for each $i=1,\dots, m$. Thus the local Milnor fibres have the homology groups of spheres.

\begin{definition}\label{vanishing cycle} Suppose that $\lambda(p_i)\leq n$. A generator
$$\gamma_i^{+}\subset \tilde{H}_{n-\lambda(p_i)}(\mathcal{F}_{i,loc}^{+})\cong \tilde{H}_{n-\lambda(p_i)}(\mathbb{S}^{n-\lambda(p_i)})$$
is called a positive vanishing cycle with respect to $F$ at the point $p_i$.\\
Suppose that $\lambda(p_i)\geq 1$. A generator
$$\gamma_i^{-}\subset \tilde{H}_{\lambda(p_i)-1}(\mathcal{F}_{i,loc}^{-})\cong \tilde{H}_{\lambda(p_i)-1}(\mathbb{S}^{\lambda(p_i)-1})$$
is called a negative vanishing cycle with respect to $F$ at the point $p_i$.\\
If $\lambda(p_i)=n+1$ one says that $F$ has no positive vanishing cycles at $p_i$. If $\lambda(p_i)=0$ one says that $F$ has no negative vanishing cycles at $p_i$.
\end{definition}

In contrast to the situation for complex hypersurface singularities this notion of vanishing cycle depends of course on the choice of morsification.

If the conditions of the Theorem \hyperref[main theorem two]{\ref*{main theorem two}} are satisfied then a set of positive, or equivalently negative, real vanishing cycles is a generating set of the top-dimensional homology groups of the real positive and negative Milnor fibres.

\begin{corollary} Suppose that the conditions of Theorem \hyperref[main theorem two]{\ref*{main theorem two}} are satisfied. Then $\gamma_1^{-},\dots, \gamma_m^{-}$ is a set of generators of the homology group 
$$\tilde{H}_k(\mathcal{F}_{\eta}^{+})\cong \tilde{H}_k(\mathcal{F}_{\eta}^{-}),\qquad k=(n-1)/2$$ 
\end{corollary}
\begin{proof} Consider the spaces 
$$\mathcal{Y}_i=f_t^{-1}(s_i-\epsilon_i, s_i+\epsilon_i)\cap \mathbb{B}_i$$
as in the proof of Theorem \hyperref[main theorem two]{\ref*{main theorem two}}. There are isomorphisms
$$\tilde{H}_k(\bar{\mathcal{F}}_{\eta}^{-})\cong H_{k+1}(\bar{\mathcal{X}}_{t,\eta}, \bar{\mathcal{F}}_{t,\eta}^{-})\cong \bigoplus_{i=1}^m H_{k+1}(Y_i, \mathcal{F}_{i, loc}^{-})$$
because $Y_i\sim \mathcal{F}_i\cup_{\partial\mathbb{D}^{\lambda}} \mathbb{D}^{\lambda}$ by the proof of Theorem \hyperref[main theorem two]{\ref*{main theorem two}}.
We claim that $\mathcal{Y}_i$ is contractible. The Morse Lemma \cite[Lemma 2.2]{Morse} gives that $\mathcal{Y}_i$ is homeomorphic to 
$$Q_i^{-1}(-\epsilon,\epsilon)=\{-\epsilon_i< Q_i(x_1,\dots, x_{n+1})< \epsilon_i\}$$
where $Q_i: U\to \mathbb{R}$ is a quadratic form of index $\lambda$ defined in a neighborhood of the origin $U\subset \mathbb{R}^{n+1}$ with $Q_i(0,\dots, 0)=0$. Then $Q_i^{-1}(0)\subset Q_i^{-1}(-\epsilon,\epsilon)$ is by \cite[Proposition 1.6]{durfee} a deformation retract. By the Local Conic Structure of algebraic sets $Q_i^{-1}(0)$ is contractible. Hence  $\mathcal{Y}_i$ is contractible. Since $\mathcal{F}_{i, loc}^{-}\neq\emptyset$ as $\lambda=(n+1)/2>0$ it follows from the long exact sequence of the pair $(\mathcal{Y}_i, \mathcal{F}_{i, loc}^{-})$ that 
$$H_{k+1}(\mathcal{Y}_i, \mathcal{F}_{i, loc}^{-})=\tilde{H}_k(\mathcal{F}_{i, loc}^{-}),\qquad k\geq 0.$$
hence 
$$\tilde{H}_k(\bar{\mathcal{F}}_{\eta}^{-})\cong \bigoplus_{i=1}^m \tilde{H}_k(\mathcal{F}_{i, loc}^{-}),\qquad k\geq 0$$
The only nontrivial homology group is in dimension $k=(n-1)/2$ so the right-hand side is generated by the negative vanishing cycles $\gamma_1^{-},\dots, \gamma_m^{-}$. Hence these generate the left-hand side as well.
\end{proof}

A set of vanishing cycles corresponding to an arbitrary $\mathbb{R}$-morsification does not generate the homology groups of the Milnor fibres (for an example see e.g Example \hyperref[example stupid]{\ref*{example stupid}}). This is to compare with the situation over the complex numbers, where the vanishing cycles determine not only the homology groups of the complex Milnor fibre but also its homotopy type.

\bibliographystyle{plain}
\bibliography{biblio.bib}

\end{document}